\documentclass{amsart}
\newtheorem{theor}{Theorem}[section]
\theoremstyle{definition} \newtheorem{defin}{Definition}[section]
\newtheorem{ex}{Example}[section]
\theoremstyle{remark} \newtheorem{rem}{Remark}[section]
\newcommand{\pn}{\par\noindent} \newcommand{\pmn}{\par\medskip\noindent}
\newcommand{\pbn}{\par\bigskip\noindent}
\newcommand{\z}{$\ast$} \newcommand{\ccc}{\circle{3}}
\newcommand{\cc}{\circle*{2}}
\begin{document}
\title{Zolotarev polynomials of degree 5, 6 and 7 with simple critical points
and their moduli spaces}
\author{Yury Kochetkov}\footnote{Department of Applied Mathematics,
Higher School of Economics, Moscow, Russia}
\date{}
\begin{abstract} A polynomial $p\in \mathbb{C}[z]$ with three finite values
is called the Zolotarev polynomial. For a class of such
polynomials with the given degree, given passport and simple
critical points we define a \emph{combinatorial moduli space}. A
combinatorial moduli space have the same essential properties as
analytical moduli space, but much easier to construct. We study
these objects for Zolotarev polynomials of degree 5, 6 and 7.
\end{abstract}

\email{yukochetkov@hse.ru, yuyukochetkov@gmail.com} \maketitle

\section{Introduction}
\pn Let $p\in\mathbb{C}[z]$ be a Zolotarev polynomial, i.e. a
polynomial with three finite pairwise different critical values
$\alpha$, $\beta$ and $\gamma$. Points $\alpha,\beta,\gamma$ are
vertices of a triangle $\Delta$ in the complex plane and the
inverse image $p^{-1}(\Delta)$ is a system of curvilinear
triangles (ovals). Any two ovals in this system are either
disjoint or have one common vertex (this common vertex is a
critical point of $p$). A simple critical point is a common vertex
of exactly two different ovals. In what follows such systems will
be called cacti [2], [3].
\begin{rem} In what follows by Galois group we understand
\emph{the geometric Galois group} [1]. \end{rem}
\begin{ex} Let
$$p=6\,\int\,(x^2+1)(x-1)(x-a)\,dx,\text{ where } a=-\frac 23+
\frac{\sqrt{35}\,i}{3}\,.$$ Then
$$p(i)\approx -6.1-0.4\,i,\,p(-i)\approx 7+5.3\,i,\,
p(1)=p(a)\approx -3+5.7\,i.$$ In the figure below the cactus $p^{-1}(\Delta)$
is presented.
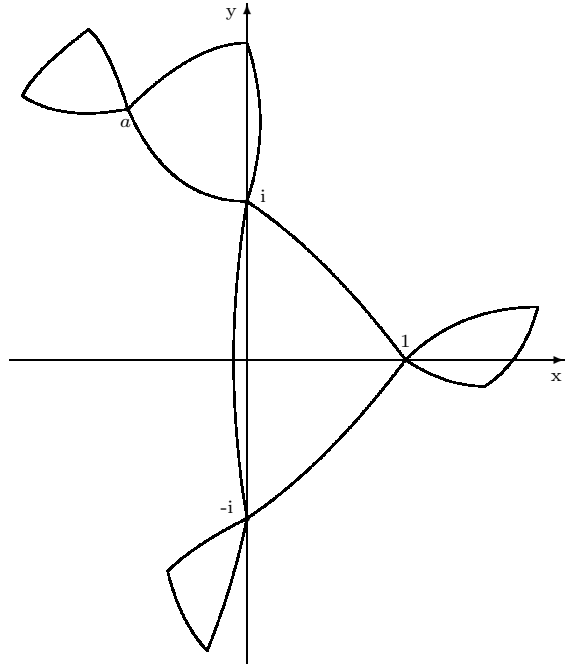
\[\begin{picture}(220,255) \put(0,120){\vector(1,0){210}}
\put(205,112){\scriptsize x} \put(90,5){\vector(0,1){250}}
\put(82,250){\scriptsize y} \put(150,120){\circle*{2}}
\put(148,125){\scriptsize 1} \put(80,62){\scriptsize -i}
\put(95,180){\scriptsize i} \put(90,60){\circle*{2}}
\put(90,180){\circle*{2}} \put(45,215){\circle*{2}}
\put(42,208){\scriptsize $a$} \linethickness{0.2mm}
\qbezier(5,220)(10,230)(30,245) \qbezier(30,245)(37,240)(45,215)
\qbezier(45,215)(20,210)(5,220) \qbezier(45,215)(70,240)(90,240)
\qbezier(90,240)(100,210)(90,180) \qbezier(45,215)(60,180)(90,180)
\qbezier(90,180)(80,120)(90,60)
\qbezier(90,180)(120,160))(150,120)
\qbezier(150,120)(120,80)(90,60)
\qbezier(150,120)(165,110)(180,110)
\qbezier(180,110)(195,120)(200,140)
\qbezier(150,120)(170,140)(200,140) \qbezier(90,60)(70,50)(60,40)
\qbezier(60,40)(65,20)(75,10) \qbezier(75,10)(85,35)(90,60)
\end{picture}\] \end{ex} \pn \begin{defin} Let $p$ be a Zolotarev polynomial
with simple critical points and let $\alpha,\beta,\gamma$ be its
critical values. If the inverse image of $\alpha$ contains $k$
critical points, the inverse image of $\beta$ contains $l$
critical points and the inverse image of $\gamma$ contains $m$
critical points, then the array $(k,l,m)$ will be called the
\emph{passport} of $p$. \end{defin} \pn Thus, the passport of the
polynomial in the above example is $(1,1,2)$. \pmn Let us note now
that the space of Zolotarev polynomials with the given passport
$(k,l,m)$ depends on one complex parameter. Indeed, let
$0,x_1,\ldots,x_{k-1}$ be critical points with the same critical
value, $1,y_1,\ldots,y_{l-1}$ be critical points with another
critical value and $z_1,\ldots,z_m$ be critical points with the
third critical value. Thus we have $k+l+m-2$ parameters. But we
have $k+l+m-3$ relations: $k-1$ relations of the type
$p(x_i)=p(0)$, $l-1$ relations of the type $p(y_i)=p(1)$ and $m-1$
relations of the type $p(z_i)=p(z_1)$. Thus, each polynomial from
this space is defined by some value of complex parameter $t$ and
the domain of $t$ is an algebraic curve $C$, maybe reducible. If
$C'$ is a connected component of $C$, then the change of $t$ in it
induces a change of the triangle $\Delta$. \pmn The generic
situation is as follows. There is a system $L$ of one-dimensional
lines in the domain $C'$. When $t\in L$, then the triangle
$\Delta$ degenerates into a line segment: some critical value
($\alpha$, for example) is strictly inside the line segment
$[\beta,\gamma]$ bounded by two other critical values. In this
case the cactus degenerates into a tree with the number of edges
equal to the double polynomial degree. Thus, we have a map in
$C'$: $L$ is the set of edges that split $C'$ into disjoint areas
--- cells. If $t_1$ and $t_2$ belong to one cell, then the
corresponding cactuses are combinatorially equivalent. The passing
of $t$ from one cell into another induces a transformation of the
cactus.
\begin{ex} We will denote vertices of $\Delta$ by three symbols
$\ast$, $\circ$ and $\bullet$ and inverse images of vertices will
be denoted by the same symbols. Let us consider the cactus
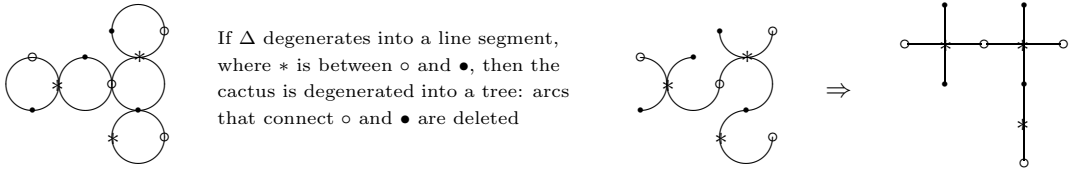
\[\begin{picture}(400,70) \multiput(10,35)(20,0){2}{\circle{20}}
\multiput(50,15)(0,20){3}{\circle{20}} \put(17,32){\z}
\put(10,45){\ccc} \put(10,25){\cc} \put(30,45){\cc}
\put(40,35){\ccc} \put(48,43){\z} \put(40,55){\cc}
\put(60,55){\ccc} \put(50,25){\cc} \put(38,12){\z}
\put(60,15){\ccc}

\put(80,50){\scriptsize If $\Delta$ degenerates into a line
segment,} \put(80,40){\scriptsize where $\ast$ is between $\circ$
and $\bullet$, then the} \put(80,30){\scriptsize cactus is
degenerated into a tree: arcs} \put(80,20){\scriptsize that
connect $\circ$ and $\bullet$ are deleted}

\put(240,35){\oval(20,20)[r]} \put(260,35){\oval(20,20)[b]}
\put(260,35){\oval(20,20)[lt]} \put(280,35){\oval(20,20)[t]}
\put(280,35){\oval(20,20)[rb]} \put(280,55){\oval(20,20)[b]}
\put(280,15){\oval(20,20)[l]} \put(280,15){\oval(20,20)[rb]}
\put(248,32){\z} \put(240,45){\ccc} \put(240,25){\cc}

\put(260,45){\cc} \put(270,35){\ccc} \put(278,43){\z}
\put(270,55){\cc} \put(290,55){\ccc} \put(280,25){\cc}
\put(268,12){\z} \put(290,15){\ccc}

\put(310,30){$\Rightarrow$}

\multiput(340,50)(30,0){3}{\ccc} \put(341,50){\line(1,0){28}}
\put(355,35){\line(0,1){30}} \put(371,50){\line(1,0){28}}
\multiput(353,47)(29,0){2}{\z} \multiput(355,35)(0,30){2}{\cc}
\put(385,6){\line(0,1){59}} \put(385,5){\ccc}
\multiput(385,35)(0,30){2}{\cc} \put(382,17){\z}
\end{picture}\] The initial cactus can be reconstructed from this tree in the
following way: to triples $\ast,\circ,\bullet$ we add arcs in such
way that the cyclic order of vertices $\ast,\circ,\bullet$ in each
thus obtained ovals is the same as in the initial cactus.
Naturally, the transformed cactus is constructed from the tree in
the same way, only the cyclic order in each constructed oval must
be opposite to the initial one.
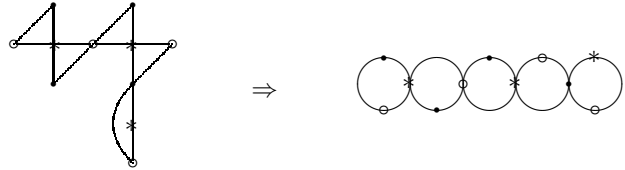
\[\begin{picture}(230,70) \multiput(0,50)(30,0){3}{\ccc}
\put(1,50){\line(1,0){28}} \put(15,35){\line(0,1){30}}
\put(31,50){\line(1,0){28}} \multiput(13,47)(29,0){2}{\z}
\multiput(15,35)(0,30){2}{\cc} \put(45,6){\line(0,1){59}}
\put(45,5){\ccc} \multiput(45,35)(0,30){2}{\cc} \put(42,17){\z}
\qbezier[25](0,50)(8,58)(15,65) \qbezier[50](15,35)(30,50)(45,65)
\qbezier[25](45,35)(53,43)(60,50) \qbezier[50](45,35)(30,20)(45,5)

\put(90,30){$\Rightarrow$}

\multiput(140,35)(20,0){5}{\circle{20}} \put(140,25){\ccc}
\put(140,45){\cc} \put(160,25){\cc} \put(170,35){\ccc}
\put(180,45){\cc} \put(200,45){\ccc} \put(210,35){\cc}
\put(220,25){\ccc} \put(147,33){\z} \put(187,33){\z}
\put(217,43){\z} \end{picture}\]
\end{ex}
\pn In the generic case the border of a cell in $C'$ consists of
three arcs, because $\Delta$ can be degenerated into the segment
$[\circ,\bullet]$, or into the segment $[\ast,\circ]$, or into the
segment $[\ast,\bullet]$.  To each arc a tree, i.e. a degenerated
cactus is corresponded. These degenerated cactuses constitute
\emph{the combinatorial border} of the initial cactus and the
corresponded transformed cactuses constitute the
\emph{combinatorial neighborhood} of the initial cactus.
\begin{ex} Let us find the combinatorial border and neighborhood of the cactus
\parbox{5cm}{\begin{picture}(100,35) \multiput(20,10)(16,0){5}{\circle{16}}
\put(52,26){\circle{16}} \put(26,8){\z} \put(42,23){\z}
\put(58,8){\z} \put(82,0){\z} \put(20,18){\ccc} \put(36,2){\ccc}
\put(52,18){\ccc} \put(68,2){\ccc} \put(84,18){\ccc}
\put(20,2){\cc} \put(44,10){\cc} \put(60,26){\cc} \put(76,10){\cc}
\end{picture}} We have:
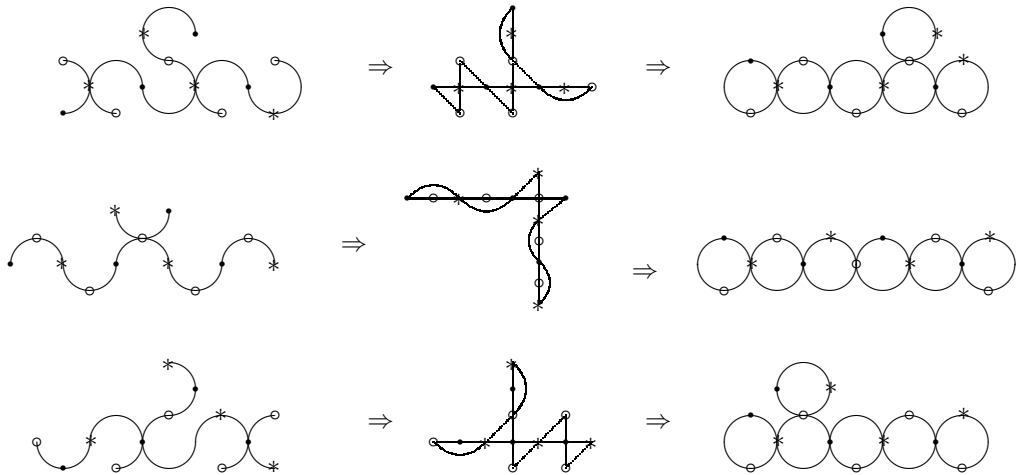
\[\begin{picture}(360,50) \put(10,15){\oval(20,20)[r]}
\put(30,15){\oval(20,20)[l]} \put(30,15){\oval(20,20)[rt]}
\put(50,15){\oval(20,20)[b]} \put(50,15){\oval(20,20)[rt]}
\put(50,35){\oval(20,20)[l]} \put(50,35){\oval(20,20)[rt]}
\put(70,15){\oval(20,20)[l]} \put(70,15){\oval(20,20)[rt]}
\put(90,15){\oval(20,20)[b]} \put(90,15){\oval(20,20)[rt]}

\put(10,25){\ccc} \put(10,5){\cc} \put(30,5){\ccc}
\put(40,15){\cc} \put(50,25){\ccc} \put(60,35){\cc}
\put(70,5){\ccc} \put(80,15){\cc} \put(90,25){\ccc}
\put(17,13){\z} \put(38,33){\z} \put(57,13){\z} \put(87,2){\z}

\put(125,20){$\Rightarrow$}

\put(150,15){\line(1,0){59}} \put(160,6){\line(0,1){18}}
\put(180,6){\line(0,1){18}} \put(180,26){\line(0,1){19}}
\put(150,15){\cc} \put(160,5){\ccc} \put(160,25){\ccc}
\put(170,15){\cc} \put(180,5){\ccc} \put(180,25){\ccc}
\put(180,45){\cc} \put(190,15){\cc} \put(210,15){\ccc}
\put(157,12){\z} \put(177,12){\z} \put(177,33){\z}
\put(197,12){\z} \qbezier[15](150,15)(155,10)(160,5)
\qbezier[30](160,25)(170,15)(180,5)
\qbezier[40](180,45)(170,35)(180,25)
\qbezier[15](180,25)(185,20)(190,15)
\qbezier[40](190,15)(200,5)(210,15)

\put(230,20){$\Rightarrow$}

\multiput(270,15)(20,0){5}{\circle{20}} \put(330,35){\circle{20}}
\put(278,13){\z} \put(318,13){\z} \put(338,33){\z}
\put(348,23){\z} \put(270,25){\cc} \put(270,5){\ccc}
\put(290,25){\ccc} \put(300,15){\cc} \put(310,5){\ccc}
\put(330,25){\ccc} \put(320,35){\cc} \put(340,15){\cc}
\put(350,5){\ccc} \end{picture}\]

\[\begin{picture}(380,50) \put(10,15){\oval(20,20)[t]}
\put(30,15){\oval(20,20)[b]} \put(50,15){\oval(20,20)[t]}
\put(50,35){\oval(20,20)[b]} \put(70,15){\oval(20,20)[b]}
\put(90,15){\oval(20,20)[t]}

\put(10,25){\ccc} \put(0,15){\cc} \put(30,5){\ccc}
\put(40,15){\cc} \put(50,25){\ccc} \put(60,35){\cc}
\put(70,5){\ccc} \put(80,15){\cc} \put(90,25){\ccc}
\put(17,13){\z} \put(37,33){\z} \put(57,13){\z} \put(97,12){\z}

\put(125,20){$\Rightarrow$}

\put(150,40){\line(1,0){60}} \put(200,0){\line(0,1){50}}
\put(150,40){\cc} \put(160,40){\ccc} \put(180,40){\ccc}
\put(190,40){\cc} \put(200,40){\ccc} \put(210,40){\cc}
\put(200,8){\ccc} \put(200,16){\cc} \put(200,24){\ccc}
\put(167,37){\z} \put(197,47){\z} \put(197,29){\z}
\put(197,-3){\z} \qbezier[40](150,40)(160,50)(170,40)
\qbezier[40](170,40)(180,30)(190,40)
\qbezier[15](190,40)(195,45)(200,50)
\qbezier[15](210,40)(205,36)(200,32)
\qbezier[35](200,32)(192,24)(200,16)
\qbezier[35](200,16)(208,8)(200,0)

\put(235,10){$\Rightarrow$}

\multiput(270,15)(20,0){6}{\circle{20}} \put(278,13){\z}
\put(308,23){\z} \put(338,13){\z} \put(368,23){\z}
\put(270,25){\cc} \put(300,15){\cc} \put(330,25){\cc}
\put(360,15){\cc} \put(270,5){\ccc} \put(290,25){\ccc}
\put(320,15){\ccc} \put(350,25){\ccc} \put(370,5){\ccc}
\end{picture}\]

\[\begin{picture}(360,50) \put(10,15){\oval(20,20)[b]}
\put(30,15){\oval(20,20)[t]} \put(30,15){\oval(20,20)[rb]}
\put(50,15){\oval(20,20)[b]} \put(50,15){\oval(20,20)[lt]}
\put(50,35){\oval(20,20)[r]} \put(70,15){\oval(20,20)[t]}
\put(70,15){\oval(20,20)[rb]} \put(90,15){\oval(20,20)[l]}

\put(0,15){\ccc} \put(10,5){\cc} \put(30,5){\ccc} \put(40,15){\cc}
\put(50,25){\ccc} \put(60,35){\cc} \put(70,5){\ccc}
\put(80,15){\cc} \put(90,25){\ccc} \put(18,13){\z} \put(47,42){\z}
\put(67,23){\z} \put(87,3){\z}

\put(125,20){$\Rightarrow$}

\put(151,15){\line(1,0){59}} \put(180,6){\line(0,1){39}}
\put(200,6){\line(0,1){18}} \put(150,15){\ccc} \put(160,15){\cc}
\put(180,15){\cc} \put(180,5){\ccc} \put(180,25){\ccc}
\put(180,35){\cc} \put(200,5){\ccc} \put(200,15){\cc}
\put(200,25){\ccc} \put(167,12){\z} \put(187,12){\z}
\put(207,12){\z} \put(177,42){\z}
\qbezier[40](150,15)(160,5)(170,15)
\qbezier[15](170,15)(175,20)(180,25)
\qbezier[40](180,25)(190,35)(180,45)
\qbezier[30](180,5)(190,15)(200,25)
\qbezier[15](200,5)(205,10)(210,15)

\put(230,20){$\Rightarrow$}

\multiput(270,15)(20,0){5}{\circle{20}} \put(290,35){\circle{20}}
\put(278,13){\z} \put(298,33){\z} \put(318,13){\z}
\put(348,23){\z} \put(270,25){\cc} \put(280,35){\cc}
\put(300,15){\cc} \put(340,15){\cc} \put(270,5){\ccc}
\put(290,25){\ccc} \put(310,5){\ccc} \put(330,25){\ccc}
\put(350,5){\ccc} \end{picture}\] \end{ex} \pn Further
degeneration transforms a line segment with a vertex inside (i.e.
a degenerated $\Delta$) into a line segment. Thus, the degenerated
triangle $\Delta$ \begin{picture}(30,4) \put(6,2){\line(1,0){19}}
\put(5,2){\ccc} \put(25,2){\cc} \put(13,-1){\z}\end{picture} can
be transformed into segment \begin{picture}(20,4)
\put(5,2){\line(1,0){10}} \put(3,-1){\z} \put(15,2){\cc}
\end{picture}, or into segment
\begin{picture}(20,4) \put(6,2){\line(1,0){9}} \put(5,2){\ccc}
\put(13,-1){\z}
\end{picture}. In the first case in a degenerated cactus we contract subsegments
\begin{picture}(20,4) \put(6,2){\line(1,0){9}} \put(5,2){\ccc}
\put(13,-1){\z} \end{picture}. And in the second case we contract
in it subsegments
\begin{picture}(20,4) \put(5,2){\line(1,0){10}} \put(15,2){\cc}
\put(3,-1){\z} \end{picture}. It means that we transforms a
3-colored tree (a degenerated cactus) into 2-colored tree. Thus,
from the degenerated cactus
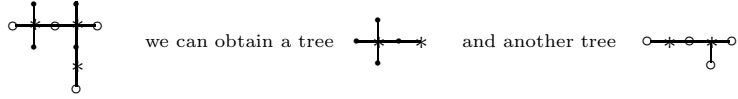
\[\begin{picture}(270,40) \put(1,28){\line(1,0){30}}
\put(8,20){\line(0,1){16}} \put(24,5){\line(0,1){32}}
\multiput(0,28)(16,0){3}{\ccc} \multiput(8,20)(0,16){2}{\cc}
\multiput(24,20)(0,16){2}{\cc} \put(24,4){\ccc}
\multiput(6,26)(16,0){2}{\z} \put(22,10){\z}

\put(50,20){\scriptsize we can obtain a tree}

\put(130,22){\line(1,0){24}} \put(138,14){\line(0,1){16}}
\multiput(130,22)(16,0){2}{\cc} \multiput(138,14)(0,16){2}{\cc}
\multiput(136,19)(16,0){2}{\z}

\put(170,20){\scriptsize and another tree}

\put(241,22){\line(1,0){30}} \put(264,14){\line(0,1){8}}
\multiput(240,22)(16,0){3}{\ccc} \put(264,13){\ccc}
\multiput(246,19)(16,0){2}{\z} \end{picture}\]
\begin{rem} A polynomial $p$ with exactly two finite critical
values is called the Shabat polynomial. If $\alpha$ and $\beta$
are its critical values, then the inverse image
$p^{-1}[\alpha,\beta]$ is a tree with number of edges equal to the
degree of $p$. Zolotarev polynomial became Shabat polynomials
exactly for those $t\in C'$ that are intersection points of two
curves from $L$. \end{rem} \pn Given the degree and the passport
of polynomials, we construct the set of cacti. For each cactus we
construct its combinatorial border
--- three degenerated cactuses, i.e. 3-colored trees with number
of edges equal to the double degree of polynomials. Two cactuses
are neighbors, if some degenerated cactus belongs to both borders.
Given a degenerated cactus we construct two 2-colored trees with
the number of edges equal to the degree of polynomials. Thus, we
have a graph: edges are enumerated by degenerated cactuses,
vertices --- by 2-colored trees, that can be obtained from
degenerated cactuses by contraction half of edges. A vertex and an
edge are incident to each other, if the corresponding 2-colored
tree can be obtained by contraction from 3-colored tree.
\begin{defin} Let $Z$ be some class of Zolotarev polynomials of the given
degree, with given passport and with simple critical points. And
let $G$ be the graph of its degenerated cactuses and 2-colored
trees. \emph{The combinatorial moduli space} of $Z$ is a map in
some surface $S$, generated by the embedding $G$ into $S$.
Moreover, faces of this map are enumerated by cactuses and
relations "being neighbors" and "being element of a border" are
preserved.
\end{defin}
\pn  Our aim is to construct combinatorial moduli spaces for some
classes of Zolotarev polynomials. Main results are as follows.
\begin{itemize}
\item A combinatorial moduli space of Zolotarev polynomials of
degree 5 and passport $(1,1,2)$ is a projective map with 5 faces,
8 edges and 4 vertices (Theorem 2.1). \item A combinatorial moduli
space of Zolotarev polynomials of degree 6, with passport
$(2,2,1)$ and the Galois group $S_6$ is a spherical map with 12
faces, 18 edges and 8 vertices (Theorem 4.1). \item A
combinatorial moduli space of Zolotarev polynomials of degree 7,
with passport $(2,2,2)$ and the Galois group $A_7$ consists of two
components. Each component is a projective map with 7 faces, 11
edges and 5 vertices (Theorem 6.1).
\end{itemize}

\section{Degree 5}
\pn Here we study Zolotarev polynomials of degree 5 and with the
passport $(1,1,2)$. The vertex of $\Delta$ that is the image of
two critical points will be denoted by star $\ast$. The vertex
following $\ast$ in the counterclockwise going around $\Delta$
will be denoted by white circle $\circ$ and the third vertex will
be denoted by black circle $\bullet$. Our cactus contains 5 ovals
and 4 points of their contact. Two of these points are inverse
images of $\ast$. In what follows inverse images of the star
vertex will be denoted by $\ast$. The same is true for inverse
images of $\circ$ and $\bullet$. \pmn We can construct 2 cacti
from 5 ovals:
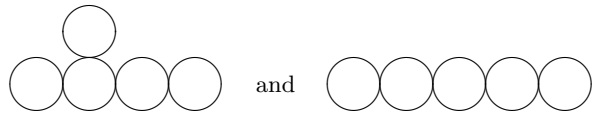
\[\begin{picture}(220,50) \multiput(10,15)(20,0){4}{\circle{20}}
\put(30,35){\circle{20}} \put(93,12){\small and}
\multiput(130,15)(20,0){5}{\circle{20}} \end{picture}\] In the
first case positions of contact points $\ast$ can be
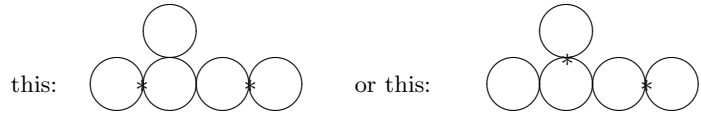
\[\begin{picture}(260,50) \put(0,12){\small this:}
\multiput(40,15)(20,0){4}{\circle{20}} \put(60,35){\circle{20}}
\put(47,12){$\ast$} \put(88,12){$\ast$}\put(130,12){\small or
this:} \multiput(190,15)(20,0){4}{\circle{20}}
\put(210,35){\circle{20}} \put(208,22){$\ast$}
\put(238,12){$\ast$} \end{picture}\] In the second case we also
have two possibilities
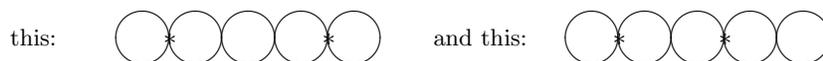
\[\begin{picture}(300,30) \put(0,12){\small this:}
\multiput(50,15)(20,0){5}{\circle{20}} \put(160,12){\small and
this:} \multiput(220,15)(20,0){5}{\circle{20}} \put(58,12){$\ast$}
\put(118,12){$\ast$} \put(228,12){$\ast$}  \put(268,12){$\ast$}
\end{picture}\] In result we have 5 cacti:
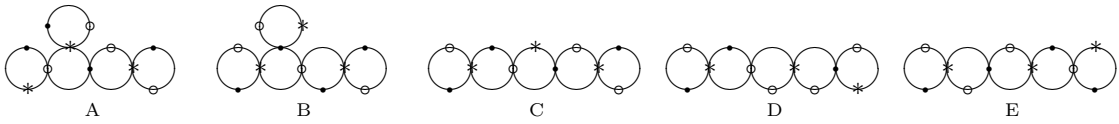
\[\begin{picture}(420,50) \multiput(8,20)(16,0){4}{\circle{16}}
\put(24,36){\circle{16}} \put(22,26){\z} \put(46,18){\z}
\put(6,10){\z} \put(16,20){\ccc} \put(8,28){\cc} \put(16,36){\cc}
\put(32,36){\ccc} \put(32,20){\cc} \put(40,28){\ccc}
\put(56,12){\ccc} \put(56,28){\cc} \put(30,2){\scriptsize A}

\multiput(88,20)(16,0){4}{\circle{16}} \put(104,36){\circle{16}}
\put(110,34){\z} \put(126,18){\z} \put(94,18){\z} \put(88,12){\cc}
\put(88,28){\ccc} \put(104,28){\cc} \put(96,36){\ccc}
\put(120,12){\cc} \put(112,20){\ccc} \put(136,12){\ccc}
\put(136,28){\cc} \put(110,2){\scriptsize B}

\multiput(168,20)(16,0){5}{\circle{16}} \put(198,26){\z}
\put(174,18){\z} \put(222,18){\z} \put(168,28){\ccc}
\put(168,12){\cc} \put(184,28){\cc} \put(192,20){\ccc}
\put(208,20){\cc} \put(216,28){\ccc} \put(232,12){\ccc}
\put(232,28){\cc} \put(198,2){\scriptsize C}

\multiput(258,20)(16,0){5}{\circle{16}} \put(264,18){\z}
\put(296,18){\z} \put(320,10){\z} \put(258,28){\ccc}
\put(258,12){\cc} \put(274,28){\cc} \put(282,20){\ccc}
\put(290,12){\ccc} \put(306,12){\ccc} \put(314,20){\cc}
\put(322,28){\ccc} \put(288,2){\scriptsize D}

\multiput(348,20)(16,0){5}{\circle{16}} \put(354,18){\z}
\put(386,18){\z} \put(410,26){\z} \put(348,28){\ccc}
\put(348,12){\cc} \put(364,12){\ccc} \put(372,20){\cc}
\put(380,28){\ccc} \put(396,28){\cc} \put(404,20){\ccc}
\put(412,12){\cc} \put(378,2){\scriptsize E}
\end{picture}\] Let us construct combinatorial borders of our
cacti. We do not mark vertices of degenerated cacti, because it is
easy to recover them.
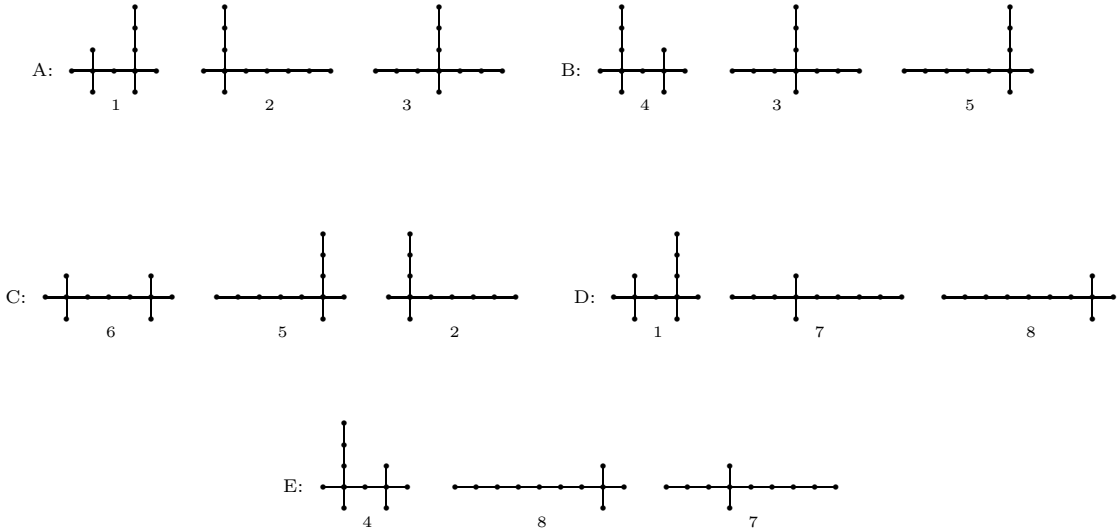
\[\begin{picture}(400,55) \put(0,18){\scriptsize A:}
\put(15,20){\line(1,0){32}} \put(23,12){\line(0,1){16}}
\put(39,12){\line(0,1){32}} \multiput(15,20)(8,0){5}{\cc}
\multiput(23,12)(0,16){2}{\cc} \multiput(39,12)(0,8){5}{\cc}
\put(30,5){\tiny 1}

\put(65,20){\line(1,0){48}} \put(73,12){\line(0,1){32}}
\multiput(65,20)(8,0){7}{\cc} \multiput(73,12)(0,8){5}{\cc}
\put(88,5){\tiny 2}

\put(130,20){\line(1,0){48}} \put(154,12){\line(0,1){32}}
\multiput(130,20)(8,0){7}{\cc} \multiput(154,12)(0,8){5}{\cc}
\put(140,5){\tiny 3}

\put(200,18){\scriptsize B:} \put(215,20){\line(1,0){32}}
\put(223,12){\line(0,1){32}} \put(239,12){\line(0,1){16}}
\multiput(215,20)(8,0){5}{\cc} \multiput(223,12)(0,8){5}{\cc}
\multiput(239,12)(0,16){2}{\cc} \put(230,5){\tiny 4}

\put(265,20){\line(1,0){48}} \put(289,12){\line(0,1){32}}
\multiput(265,20)(8,0){7}{\cc} \multiput(289,12)(0,8){5}{\cc}
\put(280,5){\tiny 3}

\put(330,20){\line(1,0){48}} \put(370,12){\line(0,1){32}}
\multiput(330,20)(8,0){7}{\cc} \multiput(370,12)(0,8){5}{\cc}
\put(353,5){\tiny 5}
\end{picture}\]

\[\begin{picture}(420,55) \put(0,18){\scriptsize C:}
\put(15,20){\line(1,0){48}} \put(23,12){\line(0,1){16}}
\put(55,12){\line(0,1){16}} \multiput(15,20)(8,0){7}{\cc}
\multiput(23,12)(0,16){2}{\cc} \multiput(55,12)(0,16){2}{\cc}
\put(38,5){\tiny 6}

\put(80,20){\line(1,0){48}} \put(120,12){\line(0,1){32}}
\multiput(80,20)(8,0){7}{\cc} \multiput(120,12)(0,8){5}{\cc}
\put(103,5){\tiny 5}

\put(145,20){\line(1,0){48}} \put(153,12){\line(0,1){32}}
\multiput(145,20)(8,0){7}{\cc} \multiput(153,12)(0,8){5}{\cc}
\put(168,5){\tiny 2}

\put(215,18){\scriptsize D:}  \put(230,20){\line(1,0){32}}
\put(238,12){\line(0,1){16}} \put(254,12){\line(0,1){32}}
\multiput(230,20)(8,0){5}{\cc} \multiput(238,12)(0,16){2}{\cc}
\multiput(254,12)(0,8){5}{\cc} \put(245,5){\tiny 1}

\put(275,20){\line(1,0){64}} \put(299,12){\line(0,1){16}}
\multiput(275,20)(8,0){9}{\cc} \multiput(299,12)(0,16){2}{\cc}
\put(306,5){\tiny 7}

\put(355,20){\line(1,0){64}} \put(411,12){\line(0,1){16}}
\multiput(355,20)(8,0){9}{\cc} \multiput(411,12)(0,16){2}{\cc}
\put(386,5){\tiny 8} \end{picture}\]

\[\begin{picture}(210,45) \put(0,18){\scriptsize E:}
\put(15,20){\line(1,0){32}} \put(23,12){\line(0,1){32}}
\put(39,12){\line(0,1){16}} \multiput(15,20)(8,0){5}{\cc}
\multiput(23,12)(0,8){5}{\cc} \multiput(39,12)(0,16){2}{\cc}
\put(30,5){\tiny 4}

\put(65,20){\line(1,0){64}} \put(121,12){\line(0,1){16}}
\multiput(65,20)(8,0){9}{\cc} \multiput(121,12)(0,16){2}{\cc}
\put(96,5){\tiny 8}

\put(145,20){\line(1,0){64}} \put(169,12){\line(0,1){16}}
\multiput(145,20)(8,0){9}{\cc} \multiput(169,12)(0,16){2}{\cc}
\put(176,5){\tiny 7} \end{picture}\] Contractions of degenerated
cacti give us four trees: \pn
$\alpha:$\parbox{15mm}{\begin{picture}(35,20)
\put(5,10){\line(1,0){24}} \put(13,2){\line(0,1){16}}
\multiput(5,10)(8,0){4}{\cc} \multiput(13,2)(0,16){2}{\cc}
\end{picture}} obtained by contraction of degenerated
cacti 1, 3, 4 and 7. \pn
$\beta:$\parbox{17mm}{\begin{picture}(40,20)
\put(4,10){\line(1,0){32}} \put(20,10){\line(0,1){8}}
\multiput(4,10)(8,0){5}{\cc} \put(20,18){\cc}
\end{picture}} obtained by contraction of degenerated
cacti 2, 3 and 5. \pn
$\gamma:$\parbox{17mm}{\begin{picture}(40,20)
\put(4,10){\line(1,0){32}} \put(12,10){\line(0,1){8}}
\multiput(4,10)(8,0){5}{\cc} \put(12,18){\cc}
\end{picture}} obtained by contraction of degenerated
cacti 1, 2, 4, 5, 6 and 8. \pn
$\delta:$\parbox{20mm}{\begin{picture}(50,20)
\put(5,10){\line(1,0){40}} \multiput(5,10)(8,0){6}{\cc}
\end{picture}} obtained by contraction of degenerated
cacti 7 and 8. \pmn Now we can formulate the statement.
\begin{theor}The combinatorial moduli space of Zolotarev polynomials
of degree 5 and with the passport $(1,1,2)$ is a projective map
with 5 faces, 8 edges and 4 vertices. \end{theor}
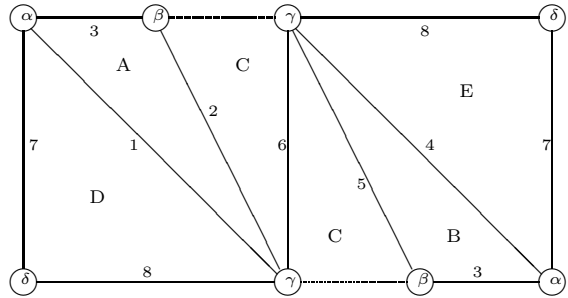
\[\begin{picture}(210,110) \multiput(5,5)(200,100){2}{\circle{10}}
\multiput(5,105)(50,0){3}{\circle{10}}
\multiput(105,5)(50,0){3}{\circle{10}}
\multiput(5,10)(100,0){3}{\line(0,1){90}}
\multiput(10,5)(100,100){2}{\line(1,0){90}}
\multiput(10,105)(150,-100){2}{\line(1,0){40}}
\multiput(8,101)(100,0){2}{\line(1,-1){93}}
\multiput(57,100)(50,0){2}{\line(1,-2){45}}
\qbezier[40](110,5)(130,5)(150,5)
\qbezier[40](60,105)(80,105)(100,105)

\put(7,55){\tiny 7} \put(201,55){\tiny 7} \put(50,7){\tiny 8}
\put(155,98){\tiny 8} \put(175,7){\tiny 3} \put(30,98){\tiny 3}
\put(45,55){\tiny 1} \put(75,68){\tiny 2} \put(101,55){\tiny 6}
\put(131,40){\tiny 5} \put(157,55){\tiny 4}

\put(30,35){\scriptsize D} \put(170,75){\scriptsize E}
\put(85,85){\scriptsize C} \put(120,20){\scriptsize C}
\put(40,85){\scriptsize A} \put(165,20){\scriptsize B}

\put(4,4){\tiny $\delta$} \put(204,104){\tiny $\delta$}
\put(4,104){\tiny $\alpha$} \put(204,4){\tiny $\alpha$}
\put(54,104){\tiny $\beta$} \put(154,4){\tiny $\beta$}
\put(104,4){\tiny $\gamma$} \put(104,104){\tiny $\gamma$}
\end{picture}\]

\section{Degree 5, analytical approach}
\pn For our polynomials (degree 5, passport $(1,1,2)$) we can easily
construct analytical moduli space. Indeed, let
$p=20\int (x^2-1)(x-b)(x-c)\,dx$. The condition $p(1)=p(-1)$ is satisfied
when $5bc=-1$. Then
$$p(1)=p(-1)=5b-\frac 1b\,,\,\,p(b)=-b^5+3b^3+2b,\,\,p(c)=\dfrac{1}{3125b^5}
-\dfrac{3}{125b^3}-\frac{2}{5b}\,.$$ The triangle $\Delta$
degenerates, if
\begin{multline*}
xy(15625x^{14}+78125x^{12}y^2+140625x^{10}y^4+78125x^8y^6-78125x^6y^8
-140625x^4y^{10}-78125x^2y^{12}-15625y^{14}-\\
-26250x^{12}-142500x^{10}y^2-333750x^8y^4-435000x^6y^6-333750x^4y^8
-142500x^2y^{10}-26250y^{12}+\\
+5775x^{10}+16125x^8y^2+10350x^6y^4-10350x^4y^6-16125x^2y^8-5775y^{10}+\\
+4660x^8+20480x^6y^2+31512x^4y^4+20480x^2y^6+4660y^8+\\
+231x^6+183x^4y^2-183x^2y^4-231y^6-42x^4-60x^2y^2-42y^4+x^2-y^2)=0
\end{multline*} where $b=x+i\,y$. The curve, defined by the above equation,
is symmetric with respect to $x$ and $y$ axes. In the figure below
the plot of this curve is presented.
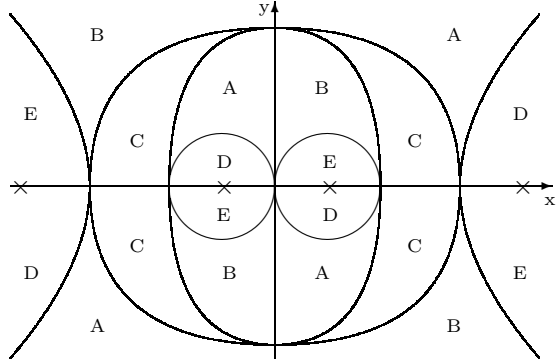
\[\begin{picture}(200,140) \put(0,70){\vector(1,0){205}}
\put(202,63){\scriptsize x} \put(100,5){\vector(0,1){135}}
\put(94,136){\scriptsize y} \put(80,70){\circle{40}}
\put(120,70){\circle{40}} \qbezier(60,70)(60,130)(100,130)
\qbezier(60,70)(60,10)(100,10) \qbezier(140,70)(140,130)(100,130)
\qbezier(140,70)(140,10)(100,10) \qbezier(30,70)(30,130)(100,130)
\qbezier(30,70)(30,10)(100,10) \qbezier(170,70)(170,130)(100,130)
\qbezier(170,70)(170,10)(100,10) \qbezier(30,70)(30,100)(0,135)
\qbezier(30,70)(30,40)(0,5) \qbezier(170,70)(170,100)(200,135)
\qbezier(170,70)(170,40)(200,5) \put(77,67){$\times$}
\put(117,67){$\times$} \put(0,67){$\times$} \put(190,67){$\times$}
\put(78,77){\scriptsize D} \put(118,77){\scriptsize E}
\put(78,57){\scriptsize E} \put(118,57){\scriptsize D}
\put(80,105){\scriptsize A} \put(115,105){\scriptsize B}
\put(80,35){\scriptsize B} \put(115,35){\scriptsize A}
\put(45,85){\scriptsize C} \put(150,85){\scriptsize C}
\put(45,45){\scriptsize C} \put(150,45){\scriptsize C}
\put(30,125){\scriptsize B} \put(165,125){\scriptsize A}
\put(30,15){\scriptsize A} \put(165,15){\scriptsize B}
\put(5,95){\scriptsize E} \put(190,95){\scriptsize D}
\put(5,35){\scriptsize D} \put(190,35){\scriptsize E}
\end{picture}\] The curve intersects $x$ axis in the origin and in
points $\pm\frac 15$ and $\pm 1$ and $y$ axis in points $\pm i/\sqrt 3$.
If $b=\pm 1,\pm \frac 15$, then $p$ is a Shabat polynomial that defines
the tree \begin{picture}(40,8) \put(4,0){\line(1,0){32}}
\put(12,0){\line(0,1){8}} \multiput(4,0)(8,0){5}{\cc} \put(12,8){\cc}
\end{picture}. If $b=\pm i/\sqrt 3$, then $p$ is a Shabat polynomial
that defines the tree \begin{picture}(40,8) \put(4,0){\line(1,0){32}}
\put(20,0){\line(0,1){8}} \multiput(4,0)(8,0){5}{\cc} \put(20,8){\cc}
\end{picture}. When $b\approx\pm 0.10$ and $b\approx\pm 1.89$
(these points are marked with crosses $\times$), then $p(b)=p(c)$
and $p$ is a Shabat polynomial that defines the tree
\begin{picture}(50,4) \put(5,2){\line(1,0){40}}
\multiput(5,2)(8,0){6}{\cc} \end{picture}. We see, that an
analytical approach gives us an "excessive" moduli space, while
combinatorial construction describes the structure of the cell
adjacency in more "compact" way.

\section{Degree 6}
\pn Here we will construct the combinatorial moduli space for
Zolotarev polynomials of degree 6 and with passport $(2,2,1)$. The
star $\ast$ will denote the vertex of the triangle $\Delta$ that
is the image of only one critical point. The vertex following
$\ast$ in the counterclockwise going around $\Delta$ will be
denoted by white circle $\circ$ and the last vertex will be
denoted by black circle $\bullet$. \pmn In all there are 15 cacti,
but three of them
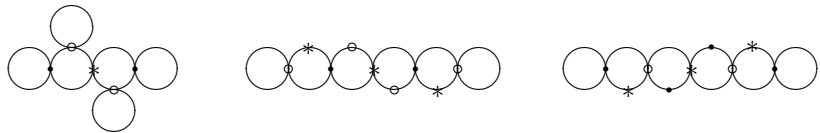
\[\begin{picture}(300,50) \multiput(8,25)(16,0){4}{\circle{16}}
\multiput(24,41)(16,-32){2}{\circle{16}} \put(30,22){\z}
\multiput(16,25)(32,0){2}{\cc} \multiput(24,33)(16,-16){2}{\ccc}

\multiput(98,25)(16,0){6}{\circle{16}} \put(136,22){\z}
\multiput(122,25)(32,0){2}{\cc} \multiput(106,25)(64,0){2}{\ccc}
\multiput(130,33)(16,-16){2}{\ccc} \put(111,30){\z}
\put(160,14){\z}

\multiput(218,25)(16,0){6}{\circle{16}} \put(256,22){\z}
\multiput(242,25)(32,0){2}{\ccc} \multiput(226,25)(64,0){2}{\cc}
\multiput(250,17)(16,16){2}{\cc} \put(232,14){\z} \put(279,31){\z}
\end{picture}\] have a nontrivial symmetry group: they admit the
rotation on $180^\circ$.\pmn We will consider "non symmetrical"
cacti with the Galois group $S_6$. There are 12 of them:
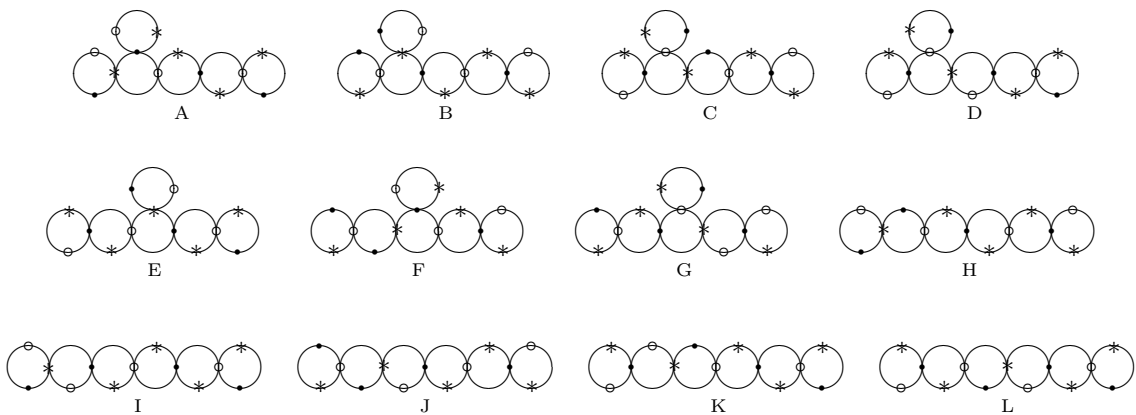
\[\begin{picture}(380,50) \multiput(8,20)(16,0){5}{\circle{16}}
\put(24,36){\circle{16}} \put(13,18){\z} \put(29,33){\z}
\put(37,25){\z} \put(53,10){\z} \put(69,25){\z} \put(8,12){\cc}
\put(24,28){\cc} \put(48,20){\cc} \put(72,12){\cc}
\put(8,28){\ccc} \put(16,36){\ccc} \put(32,20){\ccc}
\put(64,20){\ccc} \put(38,3){\scriptsize A}

\multiput(108,20)(16,0){5}{\circle{16}} \put(124,36){\circle{16}}
\put(106,10){\z} \put(122,25){\z} \put(138,10){\z}
\put(154,25){\z} \put(170,10){\z} \put(108,28){\cc}
\put(116,36){\cc} \put(132,20){\cc} \put(164,20){\cc}
\put(116,20){\ccc} \put(132,36){\ccc} \put(148,20){\ccc}
\put(172,28){\ccc} \put(138,3){\scriptsize B}

\multiput(208,20)(16,0){5}{\circle{16}} \put(224,36){\circle{16}}
\put(206,25){\z} \put(214,33){\z} \put(230,18){\z}
\put(254,25){\z} \put(270,10){\z} \put(216,20){\cc}
\put(232,36){\cc} \put(240,28){\cc} \put(264,20){\cc}
\put(208,12){\ccc} \put(224,28){\ccc} \put(248,20){\ccc}
\put(272,28){\ccc} \put(238,3){\scriptsize C}

\multiput(308,20)(16,0){5}{\circle{16}} \put(324,36){\circle{16}}
\put(306,25){\z} \put(314,34){\z} \put(330,18){\z}
\put(354,10){\z} \put(370,25){\z} \put(316,20){\cc}
\put(332,36){\cc} \put(348,20){\cc} \put(372,12){\cc}
\put(308,12){\ccc} \put(324,28){\ccc} \put(340,12){\ccc}
\put(364,20){\ccc} \put(338,3){\scriptsize D}
\end{picture}\]

\[\begin{picture}(400,40) \multiput(8,20)(16,0){5}{\circle{16}}
\put(40,36){\circle{16}} \put(6,25){\z} \put(22,10){\z}
\put(38,25){\z} \put(54,10){\z} \put(70,25){\z} \put(16,20){\cc}
\put(32,36){\cc} \put(48,20){\cc} \put(72,12){\cc}
\put(8,12){\ccc} \put(32,20){\ccc} \put(48,36){\ccc}
\put(64,20){\ccc} \put(38,3){\scriptsize E}

\multiput(108,20)(16,0){5}{\circle{16}} \put(140,36){\circle{16}}
\put(106,10){\z} \put(130,18){\z} \put(146,34){\z}
\put(154,25){\z} \put(170,10){\z} \put(108,28){\cc}
\put(124,12){\cc} \put(140,28){\cc} \put(164,20){\cc}
\put(116,20){\ccc} \put(132,36){\ccc} \put(148,20){\ccc}
\put(172,28){\ccc} \put(138,3){\scriptsize F}

\multiput(208,20)(16,0){5}{\circle{16}} \put(240,36){\circle{16}}
\put(206,10){\z} \put(222,25){\z} \put(230,34){\z}
\put(246,18){\z} \put(270,10){\z} \put(208,28){\cc}
\put(232,20){\cc} \put(248,36){\cc} \put(264,20){\cc}
\put(216,20){\ccc} \put(240,28){\ccc} \put(256,12){\ccc}
\put(272,28){\ccc} \put(238,3){\scriptsize G}

\multiput(308,20)(16,0){6}{\circle{16}} \put(314,18){\z}
\put(338,25){\z} \put(354,10){\z} \put(370,25){\z}
\put(386,10){\z} \put(308,12){\cc} \put(324,28){\cc}
\put(348,20){\cc} \put(380,20){\cc} \put(308,28){\ccc}
\put(332,20){\ccc} \put(364,20){\ccc} \put(388,28){\ccc}
\put(346,3){\scriptsize H}
\end{picture}\]

\[\begin{picture}(430,30) \multiput(8,20)(16,0){6}{\circle{16}}
\put(13,17){\z} \put(38,10){\z} \put(54,25){\z} \put(70,10){\z}
\put(86,25){\z} \put(8,12){\cc} \put(32,20){\cc} \put(64,20){\cc}
\put(88,12){\cc} \put(8,28){\ccc} \put(24,12){\ccc}
\put(48,20){\ccc} \put(80,20){\ccc} \put(48,3){\scriptsize I}

\multiput(118,20)(16,0){6}{\circle{16}} \put(116,10){\z}
\put(140,18){\z} \put(164,10){\z} \put(180,25){\z}
\put(196,10){\z} \put(118,28){\cc} \put(134,12){\cc}
\put(158,20){\cc} \put(190,20){\cc} \put(126,20){\ccc}
\put(150,12){\ccc} \put(174,20){\ccc} \put(198,28){\ccc}
\put(156,3){\scriptsize J}

\multiput(228,20)(16,0){6}{\circle{16}} \put(226,25){\z}
\put(250,18){\z} \put(274,25){\z} \put(290,10){\z}
\put(306,25){\z} \put(236,20){\cc} \put(260,28){\cc}
\put(284,20){\cc} \put(308,12){\cc} \put(228,12){\ccc}
\put(244,28){\ccc} \put(268,20){\ccc} \put(300,20){\ccc}
\put(266,3){\scriptsize K}

\multiput(338,20)(16,0){6}{\circle{16}} \put(336,25){\z}
\put(352,10){\z} \put(376,18){\z} \put(400,10){\z}
\put(416,25){\z} \put(346,20){\cc} \put(370,12){\cc}
\put(394,20){\cc} \put(418,12){\cc} \put(338,12){\ccc}
\put(362,20){\ccc} \put(386,12){\ccc} \put(410,20){\ccc}
\put(376,3){\scriptsize L} \end{picture}\] Also we have 18
degenerated cacti:
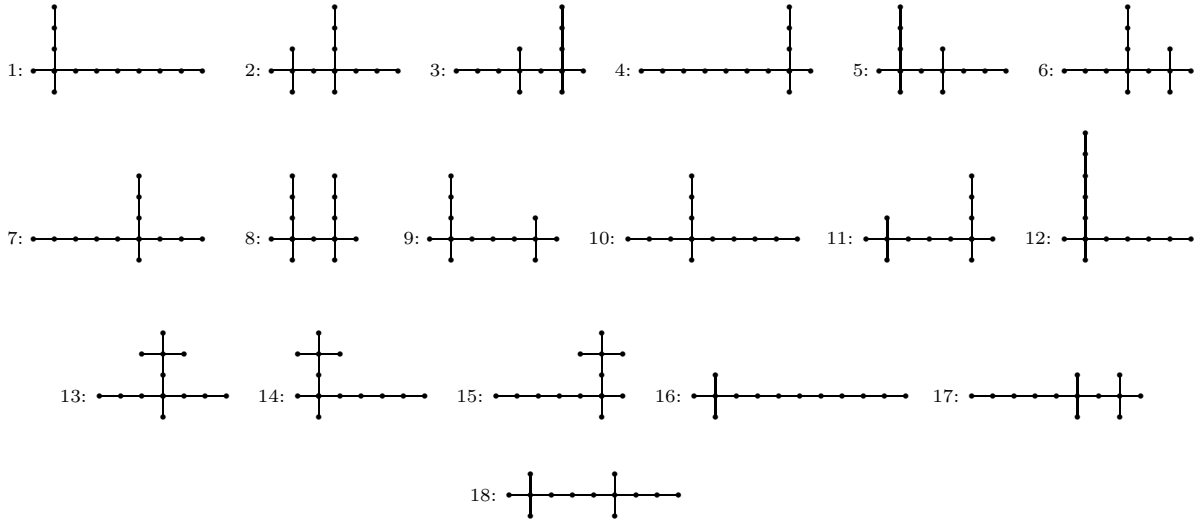
\[\begin{picture}(450,40) \put(0,8){\scriptsize 1:}
\put(10,10){\line(1,0){64}} \put(18,2){\line(0,1){32}}
\multiput(10,10)(8,0){9}{\cc} \multiput(18,2)(0,8){5}{\cc}

\put(90,8){\scriptsize 2:} \put(100,10){\line(1,0){48}}
\put(108,2){\line(0,1){16}} \put(124,2){\line(0,1){32}}
\multiput(100,10)(8,0){7}{\cc} \multiput(108,2)(0,16){2}{\cc}
\multiput(124,2)(0,8){5}{\cc}

\put(160,8){\scriptsize 3:} \put(170,10){\line(1,0){48}}
\put(194,2){\line(0,1){16}} \put(210,2){\line(0,1){32}}
\multiput(170,10)(8,0){7}{\cc} \multiput(194,2)(0,16){2}{\cc}
\multiput(210,2)(0,8){5}{\cc}

\put(230,8){\scriptsize 4:} \put(240,10){\line(1,0){64}}
\put(296,2){\line(0,1){32}} \multiput(240,10)(8,0){9}{\cc}
\multiput(296,2)(0,8){5}{\cc}

\put(320,8){\scriptsize 5:} \put(330,10){\line(1,0){48}}
\put(338,2){\line(0,1){32}} \put(354,2){\line(0,1){16}}
\multiput(330,10)(8,0){7}{\cc} \multiput(338,2)(0,8){5}{\cc}
\multiput(354,2)(0,16){2}{\cc}

\put(390,8){\scriptsize 6:} \put(400,10){\line(1,0){48}}
\put(424,2){\line(0,1){32}} \put(440,2){\line(0,1){16}}
\multiput(400,10)(8,0){7}{\cc} \multiput(424,2)(0,8){5}{\cc}
\multiput(440,2)(0,16){2}{\cc} \end{picture}\]

\[\begin{picture}(450,40) \put(0,8){\scriptsize 7:}
\put(10,10){\line(1,0){64}} \put(50,2){\line(0,1){32}}
\multiput(10,10)(8,0){9}{\cc} \multiput(50,2)(0,8){5}{\cc}

\put(90,8){\scriptsize 8:} \put(100,10){\line(1,0){32}}
\put(108,2){\line(0,1){32}} \put(124,2){\line(0,1){32}}
\multiput(100,10)(8,0){5}{\cc} \multiput(108,2)(0,8){5}{\cc}
\multiput(124,2)(0,8){5}{\cc}

\put(150,8){\scriptsize 9:} \put(160,10){\line(1,0){48}}
\put(168,2){\line(0,1){32}} \put(200,2){\line(0,1){16}}
\multiput(160,10)(8,0){7}{\cc} \multiput(168,2)(0,8){5}{\cc}
\multiput(200,2)(0,16){2}{\cc}

\put(220,8){\scriptsize 10:} \put(235,10){\line(1,0){64}}
\put(259,2){\line(0,1){32}} \multiput(235,10)(8,0){9}{\cc}
\multiput(259,2)(0,8){5}{\cc}

\put(310,8){\scriptsize 11:} \put(325,10){\line(1,0){48}}
\put(333,2){\line(0,1){16}} \put(365,2){\line(0,1){32}}
\multiput(325,10)(8,0){7}{\cc} \multiput(333,2)(0,16){2}{\cc}
\multiput(365,2)(0,8){5}{\cc}

\put(385,8){\scriptsize 12:} \put(400,10){\line(1,0){48}}
\put(408,2){\line(0,1){48}} \multiput(400,10)(8,0){7}{\cc}
\multiput(408,2)(0,8){7}{\cc}
\end{picture}\]

\[\begin{picture}(410,40) \put(0,8){\scriptsize 13:}
\put(15,10){\line(1,0){48}} \put(39,2){\line(0,1){32}}
\put(31,26){\line(1,0){16}} \multiput(15,10)(8,0){7}{\cc}
\multiput(39,2)(0,8){5}{\cc} \multiput(31,26)(16,0){2}{\cc}

\put(75,8){\scriptsize 14:} \put(90,10){\line(1,0){48}}
\put(98,2){\line(0,1){32}} \put(90,26){\line(1,0){16}}
\multiput(90,10)(8,0){7}{\cc} \multiput(98,2)(0,8){5}{\cc}
\multiput(90,26)(16,0){2}{\cc}

\put(150,8){\scriptsize 15:} \put(165,10){\line(1,0){48}}
\put(205,2){\line(0,1){32}} \put(197,26){\line(1,0){16}}
\multiput(165,10)(8,0){7}{\cc} \multiput(205,2)(0,8){5}{\cc}
\multiput(197,26)(16,0){2}{\cc}

\put(225,8){\scriptsize 16:} \put(240,10){\line(1,0){80}}
\put(248,2){\line(0,1){16}} \multiput(240,10)(8,0){11}{\cc}
\multiput(248,2)(0,16){2}{\cc}

\put(330,8){\scriptsize 17:} \put(345,10){\line(1,0){64}}
\put(385,2){\line(0,1){16}} \put(401,2){\line(0,1){16}}
\multiput(345,10)(8,0){9}{\cc} \multiput(385,2)(0,16){2}{\cc}
\multiput(401,2)(0,16){2}{\cc} \end{picture}\]

\[\begin{picture}(100,20) \put(0,8){\scriptsize 18:}
\put(15,10){\line(1,0){64}} \put(23,2){\line(0,1){16}}
\put(55,2){\line(0,1){16}} \multiput(15,10)(8,0){9}{\cc}
\multiput(23,2)(0,16){2}{\cc} \multiput(55,2)(0,16){2}{\cc}
\end{picture}\] And at last we have 8 trees --- contractions of
degenerated cacti:
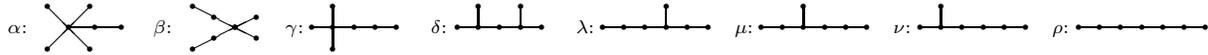
\[\begin{picture}(460,20) \put(0,8){\scriptsize $\alpha$:}
\put(15,2){\line(1,1){16}} \put(15,18){\line(1,-1){16}}
\put(23,10){\line(1,0){20}} \multiput(15,2)(16,16){2}{\cc}
\multiput(15,18)(16,-16){2}{\cc} \multiput(23,10)(10,0){3}{\cc}

\put(55,8){\scriptsize $\beta$:} \put(70,2){\line(2,1){24}}
\put(70,18){\line(2,-1){24}} \multiput(70,2)(8,4){4}{\cc}
\multiput(70,18)(8,-4){4}{\cc}

\put(105,8){\scriptsize $\gamma$:} \put(115,10){\line(1,0){32}}
\put(123,2){\line(0,1){16}} \multiput(115,10)(8,0){5}{\cc}
\multiput(123,2)(0,16){2}{\cc}

\put(160,8){\scriptsize $\delta$:} \put(170,10){\line(1,0){32}}
\multiput(178,10)(16,0){2}{\line(0,1){8}}
\multiput(170,10)(8,0){5}{\cc} \multiput(178,18)(16,0){2}{\cc}

\put(215,8){\scriptsize $\lambda$:} \put(225,10){\line(1,0){40}}
\put(249,10){\line(0,1){8}} \multiput(225,10)(8,0){6}{\cc}
\put(249,18){\cc}

\put(275,8){\scriptsize $\mu$:} \put(285,10){\line(1,0){40}}
\put(301,10){\line(0,1){8}} \multiput(285,10)(8,0){6}{\cc}
\put(301,18){\cc}

\put(335,8){\scriptsize $\nu$:} \put(345,10){\line(1,0){40}}
\put(353,10){\line(0,1){8}} \multiput(345,10)(8,0){6}{\cc}
\put(353,18){\cc}

\put(395,8){\scriptsize $\rho$:} \put(405,10){\line(1,0){48}}
\multiput(405,10)(8,0){7}{\cc} \end{picture}\] For each
degenerated cactus we: a) enumerate cacti to which border it
belongs; b) enumerate trees that can be obtained from it by
contraction.
$$\begin{array}{llllll} \hspace{2mm}1:\,\, A,J;\, \lambda,\nu.&
\hspace{2mm}2:\,\, A,E;\, \alpha,\lambda; & \hspace{2mm}3:\,\,
A,H;\, \alpha,\nu. & \hspace{2mm}4:\,\, B,K;\,\mu,\nu.&
\hspace{2mm} 5:\,\,
B,I;\,\alpha,\nu. & \hspace{2mm}6:\,B,E;\,\alpha,\mu. \\
\hspace{2mm}7:\,\,C,G;\, \beta,\lambda. &
\hspace{2mm}8:\,\,C,D;\,\beta,\delta. &
\hspace{2mm}9:\,\,C,L;\,\delta, \lambda. & 10:\,D,F;\,\beta,\mu. &
11:\, D,L;\, \delta,\mu. & 12:\,E,L;\, \lambda,\mu.\\ 13:\, F,G;\,
\beta,\gamma. & 14:\,F,K;\,\gamma,\mu. &
15:\,G,J;\,\gamma,\lambda. & 16:\, H,I;\, \nu,\rho. & 17:\,
H,I;\,\alpha,\rho. & 18:\, J,K;\,\gamma,\nu. \end{array}$$ As we
have 12 cacti, 18 degenerated cacti and 8 trees, then our map is a
spherical map. \begin{theor} The combinatorial moduli space of
Zolotarev polynomials of degree 6, with passport $(2,2,1)$ and the
Galois group $S_6$ is a spherical map presented in figure below
(here $A$ is the external face): \end{theor}
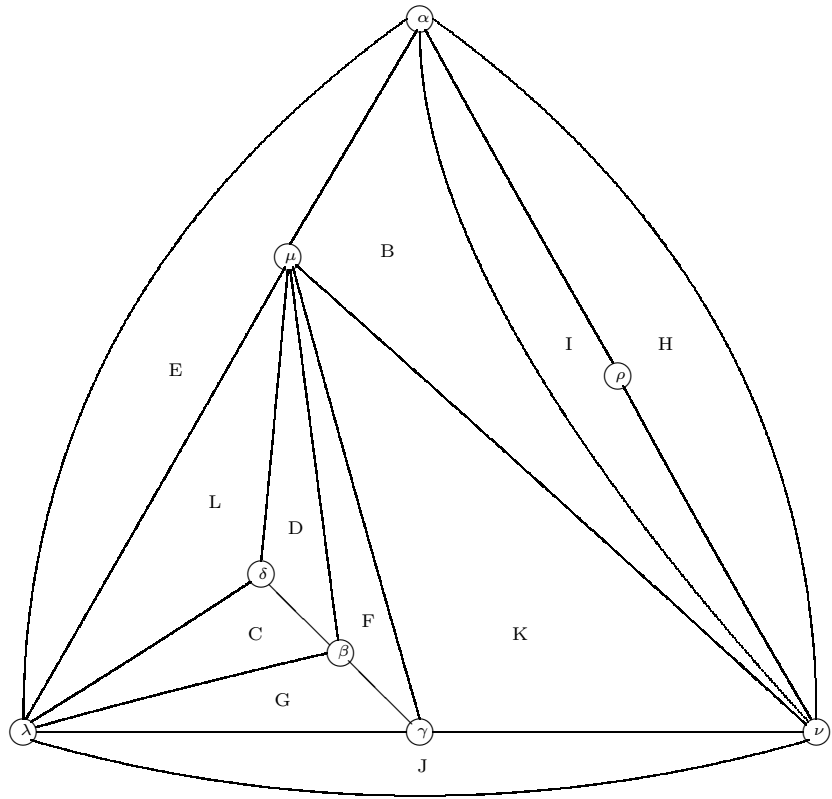
\[\begin{picture}(310,310) \multiput(5,35)(150,0){3}{\circle{10}}
\multiput(10,35)(150,0){2}{\line(1,0){140}}
\put(155,305){\circle{10}} \put(230,170){\circle{10}}
\put(105,215){\circle{10}} \qbezier(157,301)(192,237)(228,175)
\qbezier(232,166)(267,101)(303,40)
\qbezier(108,212)(206,125)(301,38) \qbezier(6,40)(55,125)(104,211)
\qbezier(106,220)(130,260)(154,301)
\qbezier(155,300)(155,200)(302,39)
\multiput(125,65)(-30,30){2}{\circle{10}}
\qbezier(105,210)(100,152)(95,100)
\qbezier(106,210)(115,140)(124,70)
\qbezier(107,211)(131,123)(155,40) \qbezier(8,39)(50,65)(91,92)
\put(152,38){\line(-1,1){24}} \put(121,68){\line(-1,1){23}}
\qbezier(10,37)(60,51)(120,65) \qbezier(5,40)(5,200)(150,305)
\qbezier(305,40)(305,200)(160,305) \qbezier(8,32)(155,-10)(302,32)

\put(4,34){\tiny $\lambda$} \put(154,34){\tiny $\gamma$}
\put(304,34){\tiny $\nu$} \put(124,64){\tiny $\beta$}
\put(94,93){\tiny $\delta$} \put(104,214){\tiny $\mu$}
\put(154,304){\tiny $\alpha$} \put(229,169){\tiny $\rho$}

\put(154,20){\scriptsize J} \put(100,45){\scriptsize G}
\put(190,70){\scriptsize K} \put(90,70){\scriptsize C}
\put(75,120){\scriptsize L} \put(105,110){\scriptsize D}
\put(133,75){\scriptsize F} \put(140,215){\scriptsize B}
\put(210,180){\scriptsize I} \put(60,170){\scriptsize E}
\put(245,180){\scriptsize H}  \end{picture}\]

\section{Degree 6, analytical approach}
\pn Let $p=\int (x^2-1)(x^2+ax+b)(x-c)\,dx$. If $p(1)=p(-1)$ and values
of $p$ in roots of $x^2+ax+b$ coincide, then $p$ is a Zolotarev polynomial
with the passport $(2,2,1)$. The set of solutions of the system
$$\left\{\begin{array}{l} p(1)=p(-1)\\ \text{the remainder $r$, where
$p=s\cdot(x^2+ax+b)+r$ has degree 0}\end{array}\right.$$ has a rational
parametrization:
$$a=\frac{\sqrt 3\,(z^2+1)}{2z}\,,\quad b=\frac{z^2-1}{2}\,,\quad
c=\frac{\sqrt 3\,(z^2+1)}{z(5z^2-3)}\,.$$ In the figure below curves
of the system $L$ in the first quadrant are presented. Axes $x$ and $y$
belong to $L$. If a segment at axis $x$ (or $y$) separates two cells, then
corresponding cactuses are mirror symmetric.
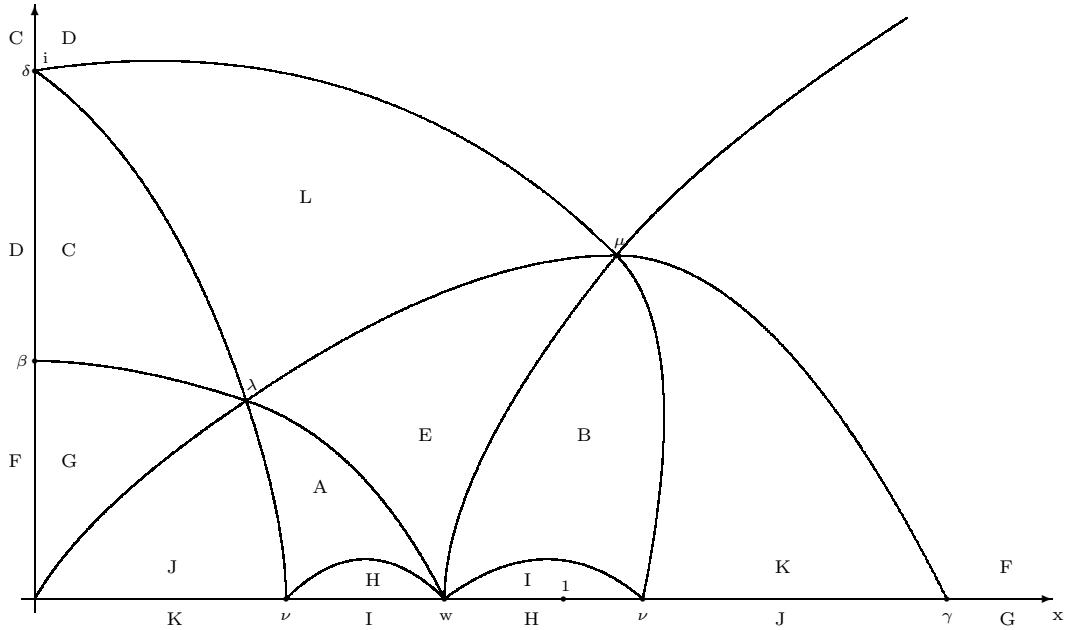
\[\begin{picture}(410,250) \put(15,20){\vector(1,0){390}}
\put(20,15){\vector(0,1){230}} \put(220,20){\circle*{2}}
\put(219,23){\tiny 1} \put(20,220){\circle*{2}} \put(23,223){\tiny
i} \put(20,110){\circle*{2}} \put(13,108){\tiny $\beta$}
\put(365,20){\circle*{2}} \put(363,12){\tiny $\gamma$}
\put(100,95){\circle*{2}} \put(115,20){\circle*{2}}
\qbezier(100,95)(115,50)(115,20) \qbezier(100,95)(70,185)(20,220)
\put(175,20){\circle*{2}} \qbezier(20,110)(55,110)(100,95)
\qbezier(100,95)(145,80)(175,20) \qbezier(115,20)(145,50)(175,20)
\put(250,20){\circle*{2}} \qbezier(175,20)(215,50)(250,20)
\put(240,150){\circle*{2}} \qbezier(20,220)(150,240)(240,150)
\qbezier(240,150)(270,120)(250,20) \qbezier(20,20)(40,55)(100,95)
\qbezier(100,95)(180,150)(240,150)
\qbezier(240,150)(300,150)(365,20)
\qbezier(175,20)(175,70)(240,150)
\qbezier(240,150)(273,190)(350,240) \put(15,218){\tiny $\delta$}
\put(100,99){\tiny $\lambda$} \put(239,154){\tiny $\mu$}
\put(113,12){\tiny $\nu$} \put(248,12){\tiny $\nu$}
\put(173,12){\tiny w} \put(10,70){\scriptsize F}
\put(30,70){\scriptsize G} \put(10,150){\scriptsize D}
\put(30,150){\scriptsize C} \put(10,230){\scriptsize C}
\put(30,230){\scriptsize D} \put(70,10){\scriptsize K}
\put(70,30){\scriptsize J} \put(145,10){\scriptsize I}
\put(145,25){\scriptsize H} \put(205,10){\scriptsize H}
\put(205,25){\scriptsize I} \put(300,10){\scriptsize J}
\put(300,30){\scriptsize K} \put(385,10){\scriptsize G}
\put(385,30){\scriptsize F} \put(405,12){\scriptsize x}
\put(225,80){\scriptsize B} \put(165,80){\scriptsize E}
\put(125,60){\scriptsize A}  \put(120,170){\scriptsize L}
\end{picture}\] Vertices have the following coordinates: $\beta
\approx 0.44\,i$; $\gamma \approx 1.73$; $\delta=i$; $\lambda
\approx 0.39+0.38\,i$; $\mu \approx 1.1+0.65\,i$; $\nu \approx
0.47 \text{ and } \approx 1.14$. At point $w=\sqrt{3/5}$ $c$ becames infinite.

\section{Degree 7}
\pn Here we will study Zolotarev polynomials of degree 7 and with
passport $(2,2,2)$. As above, vertices of $\Delta$ in the
counterclockwise going around $\Delta$ are in the order
$\ast,\circ,\bullet$. We will mark by $\ast,\circ,\bullet$ only
contact points of cacti. We have 20 of them:
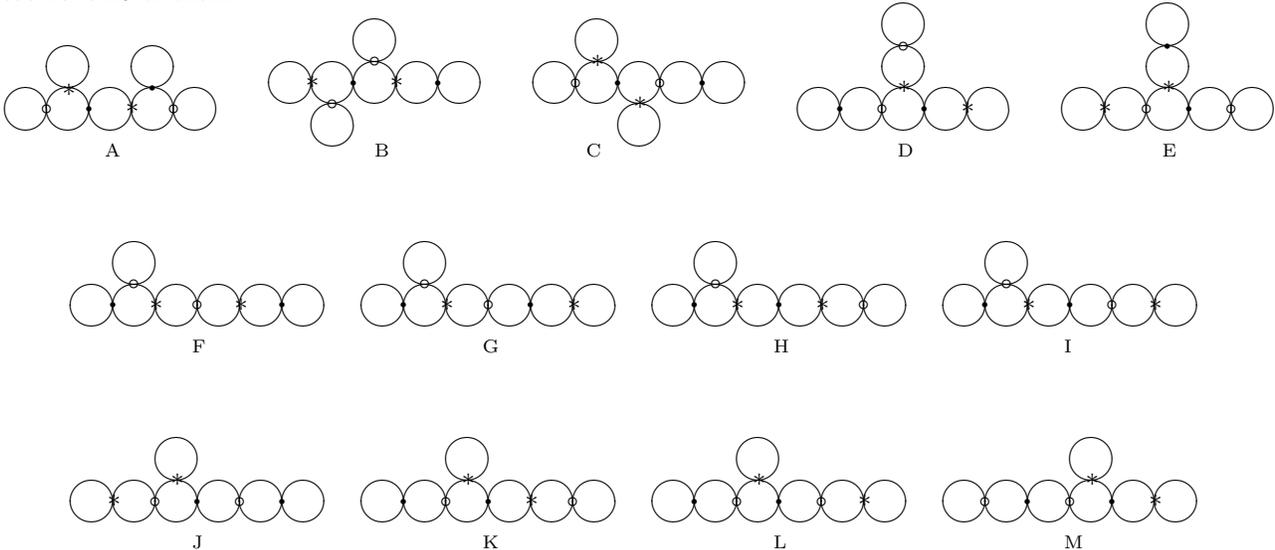
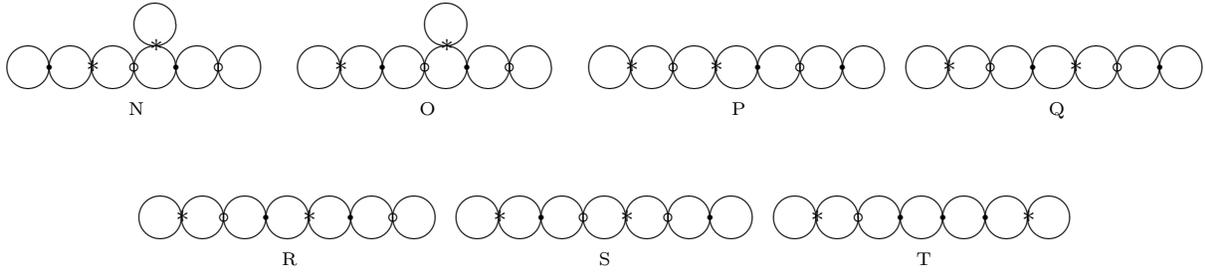
\[\begin{picture}(480,55) \multiput(8,20)(16,0){5}{\circle{16}}
\multiput(24,36)(32,0){2}{\circle{16}} \put(16,20){\ccc}
\put(22,25){\z} \put(32,20){\cc} \put(46,18){\z} \put(64,20){\ccc}
\put(56,28){\cc} \put(38,2){\scriptsize A}

\multiput(108,30)(16,0){5}{\circle{16}} \put(124,14){\circle{16}}
\put(140,46){\circle{16}} \put(114,28){\z} \put(124,22){\ccc}
\put(132,30){\cc} \put(140,38){\ccc} \put(146,28){\z}
\put(164,30){\cc} \put(140,2){\scriptsize B}

\multiput(208,30)(16,0){5}{\circle{16}} \put(224,46){\circle{16}}
\put(240,14){\circle{16}}  \put(222,36){\z} \put(216,30){\ccc}
\put(232,30){\cc} \put(238,20){\z} \put(248,30){\ccc}
\put(264,30){\cc}  \put(220,2){\scriptsize C}

\multiput(308,20)(16,0){5}{\circle{16}}
\multiput(340,36)(0,16){2}{\circle{16}} \put(338,26){\z}
\put(340,44){\ccc} \put(332,20){\ccc} \put(316,20){\cc}
\put(348,20){\cc} \put(362,18){\z} \put(338,2){\scriptsize D}

\multiput(408,20)(16,0){5}{\circle{16}}
\multiput(440,36)(0,16){2}{\circle{16}} \put(438,26){\z}
\put(440,44){\cc} \put(432,20){\ccc} \put(414,18){\z}
\put(448,20){\cc} \put(464,20){\ccc} \put(438,2){\scriptsize E}
\end{picture}\]

\[\begin{picture}(430,50) \multiput(8,20)(16,0){6}{\circle{16}}
\put(24,36){\circle{16}} \put(16,20){\cc} \put(24,28){\ccc}
\put(30,18){\z} \put(48,20){\ccc} \put(62,18){\z} \put(80,20){\cc}
\put(46,2){\scriptsize F}

\multiput(118,20)(16,0){6}{\circle{16}} \put(134,36){\circle{16}}
\put(126,20){\cc} \put(134,28){\ccc} \put(140,18){\z}
\put(158,20){\ccc} \put(188,18){\z} \put(174,20){\cc}
\put(156,2){\scriptsize G}

\multiput(228,20)(16,0){6}{\circle{16}} \put(244,36){\circle{16}}
\put(236,20){\cc} \put(244,28){\ccc} \put(250,18){\z}
\put(268,20){\cc} \put(282,18){\z} \put(300,20){\ccc}
\put(266,2){\scriptsize H}

\multiput(338,20)(16,0){6}{\circle{16}} \put(354,36){\circle{16}}
\put(346,20){\cc} \put(354,28){\ccc} \put(360,18){\z}
\put(378,20){\cc} \put(394,20){\ccc} \put(408,18){\z}
\put(376,2){\scriptsize I}
\end{picture}\]

\[\begin{picture}(430,50) \multiput(8,20)(16,0){6}{\circle{16}}
\put(40,36){\circle{16}} \put(14,18){\z} \put(32,20){\ccc}
\put(38,26){\z} \put(48,20){\cc} \put(64,20){\ccc}
\put(80,20){\cc} \put(46,2){\scriptsize J}

\multiput(118,20)(16,0){6}{\circle{16}} \put(150,36){\circle{16}}
\put(126,20){\cc} \put(142,20){\ccc} \put(148,26){\z}
\put(158,20){\cc} \put(172,18){\z} \put(190,20){\ccc}
\put(156,2){\scriptsize K}

\multiput(228,20)(16,0){6}{\circle{16}} \put(260,36){\circle{16}}
\put(236,20){\cc} \put(252,20){\ccc} \put(258,26){\z}
\put(268,20){\cc} \put(284,20){\ccc} \put(298,18){\z}
\put(266,2){\scriptsize L}

\multiput(338,20)(16,0){6}{\circle{16}} \put(386,36){\circle{16}}
\put(346,20){\ccc} \put(362,20){\cc} \put(378,20){\ccc}
\put(384,26){\z} \put(394,20){\cc} \put(408,18){\z}
\put(376,2){\scriptsize M}
\end{picture}\]

\[\begin{picture}(450,50) \multiput(8,20)(16,0){6}{\circle{16}}
\put(56,36){\circle{16}} \put(16,20){\cc} \put(30,18){\z}
\put(48,20){\ccc} \put(54,26){\z} \put(64,20){\cc}
\put(80,20){\ccc} \put(46,2){\scriptsize N}

\multiput(118,20)(16,0){6}{\circle{16}} \put(166,36){\circle{16}}
\put(124,18){\z} \put(142,20){\cc} \put(158,20){\ccc}
\put(164,26){\z} \put(174,20){\cc} \put(190,20){\ccc}
\put(156,2){\scriptsize O}

\multiput(228,20)(16,0){7}{\circle{16}} \put(234,18){\z}
\put(252,20){\ccc} \put(266,18){\z} \put(284,20){\cc}
\put(300,20){\ccc} \put(316,20){\cc} \put(274,2){\scriptsize P}

\multiput(348,20)(16,0){7}{\circle{16}} \put(354,18){\z}
\put(372,20){\ccc} \put(388,20){\cc} \put(402,18){\z}
\put(420,20){\ccc} \put(436,20){\cc} \put(394,2){\scriptsize Q}
\end{picture}\]

\[\begin{picture}(350,30) \multiput(8,20)(16,0){7}{\circle{16}}
\put(14,18){\z} \put(32,20){\ccc} \put(48,20){\cc} \put(62,18){\z}
\put(80,20){\cc} \put(96,20){\ccc} \put(54,2){\scriptsize R}

\multiput(128,20)(16,0){7}{\circle{16}} \put(134,18){\z}
\put(152,20){\cc} \put(168,20){\ccc} \put(182,18){\z}
\put(200,20){\ccc} \put(216,20){\cc} \put(174,2){\scriptsize S}

\multiput(248,20)(16,0){7}{\circle{16}} \put(254,18){\z}
\put(272,20){\ccc} \put(288,20){\cc} \put(304,20){\cc}
\put(320,20){\cc} \put(334,18){\z} \put(294,2){\scriptsize T}
\end{picture}\] The Galois group of $D,E,G,I,K,N$ is $PSL_2(7)$
of order 168. Other 14 cacti have group $A_7$ as the Galois group.
We will construct the combinatorial moduli space exactly for these
14 cacti. \pmn Now let us enumerate degenerated cacti:
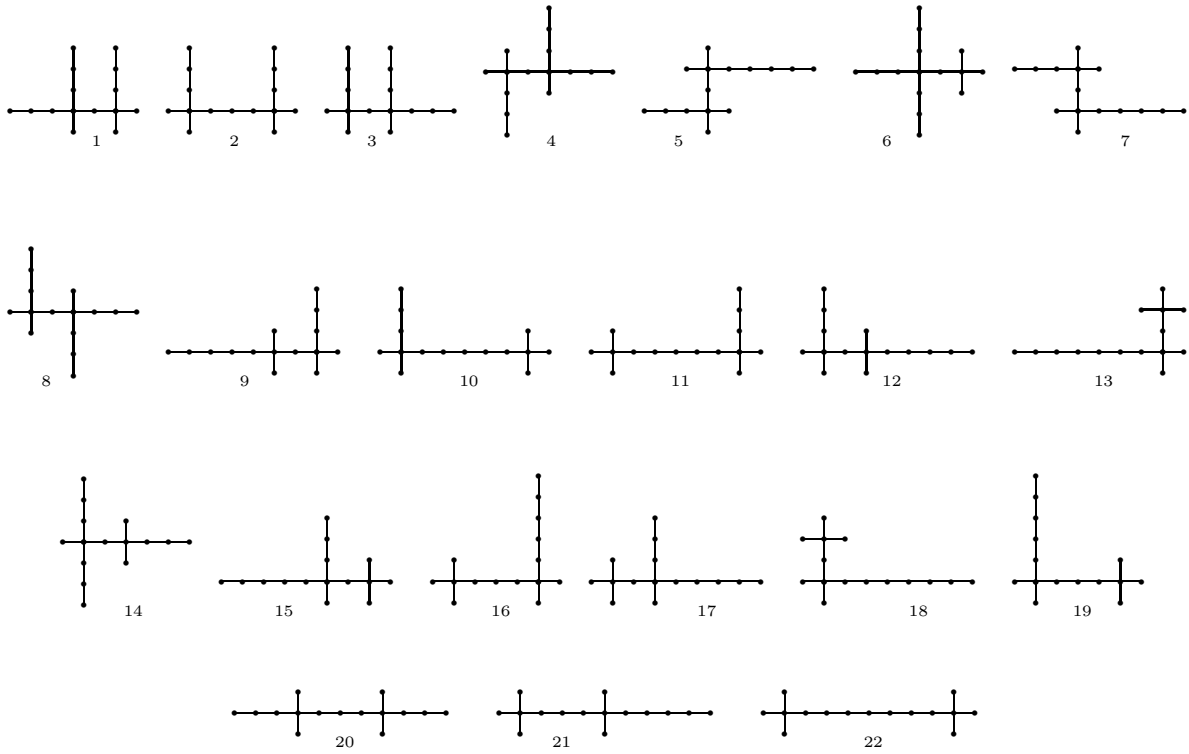
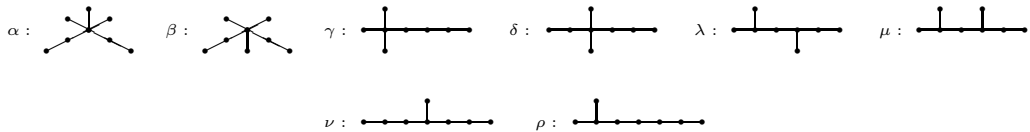
\[\begin{picture}(450,60) \put(0,15){\line(1,0){48}}
\multiput(24,7)(16,0){2}{\line(0,1){32}}
\multiput(0,15)(8,0){7}{\cc} \multiput(24,7)(0,8){5}{\cc}
\multiput(40,7)(0,8){5}{\cc} \put(31,2){\tiny 1}

\put(60,15){\line(1,0){48}}
\multiput(68,7)(32,0){2}{\line(0,1){32}}
\multiput(60,15)(8,0){7}{\cc} \multiput(68,7)(0,8){5}{\cc}
\multiput(100,7)(0,8){5}{\cc} \put(83,2){\tiny 2}

\put(120,15){\line(1,0){48}}
\multiput(128,7)(16,0){2}{\line(0,1){32}}
\multiput(120,15)(8,0){7}{\cc} \multiput(128,7)(0,8){5}{\cc}
\multiput(144,7)(0,8){5}{\cc} \put(135,2){\tiny 3}

\put(180,30){\line(1,0){48}}
\multiput(188,6)(16,16){2}{\line(0,1){32}}
\multiput(180,30)(8,0){7}{\cc} \multiput(188,6)(0,8){5}{\cc}
\multiput(204,22)(0,8){5}{\cc} \put(203,2){\tiny 4}

\put(240,15){\line(1,0){32}} \put(264,7){\line(0,1){32}}
\put(256,31){\line(1,0){48}} \multiput(240,15)(8,0){5}{\cc}
\multiput(264,7)(0,8){5}{\cc} \multiput(256,31)(8,0){7}{\cc}
\put(251,2){\tiny 5}

\put(320,30){\line(1,0){48}} \put(344,6){\line(0,1){48}}
\put(360,22){\line(0,1){16}} \multiput(320,30)(8,0){7}{\cc}
\multiput(344,6)(0,8){7}{\cc} \multiput(360,22)(0,16){2}{\cc}
\put(330,2){\tiny 6}

\put(380,31){\line(1,0){32}} \put(404,7){\line(0,1){32}}
\put(396,15){\line(1,0){48}}  \multiput(380,31)(8,0){5}{\cc}
\multiput(404,7)(0,8){5}{\cc} \multiput(396,15)(8,0){7}{\cc}
\put(420,2){\tiny 7}
\end{picture}\]

\[\begin{picture}(450,60) \put(0,30){\line(1,0){48}}
\put(8,22){\line(0,1){32}} \put(24,6){\line(0,1){32}}
\multiput(0,30)(8,0){7}{\cc} \multiput(8,22)(0,8){5}{\cc}
\multiput(24,6)(0,8){5}{\cc} \put(12,2){\tiny 8}

\put(60,15){\line(1,0){64}} \put(116,7){\line(0,1){32}}
\put(100,7){\line(0,1){16}} \multiput(60,15)(8,0){9}{\cc}
\multiput(116,7)(0,8){5}{\cc} \multiput(100,7)(0,16){2}{\cc}
\put(87,2){\tiny 9}

\put(140,15){\line(1,0){64}} \put(148,7){\line(0,1){32}}
\put(196,7){\line(0,1){16}} \multiput(140,15)(8,0){9}{\cc}
\multiput(148,7)(0,8){5}{\cc} \multiput(196,7)(0,16){2}{\cc}
\put(170,2){\tiny 10}

\put(220,15){\line(1,0){64}} \put(228,7){\line(0,1){16}}
\put(276,7){\line(0,1){32}} \multiput(220,15)(8,0){9}{\cc}
\multiput(228,7)(0,16){2}{\cc} \multiput(276,7)(0,8){5}{\cc}
\put(250,2){\tiny 11}

\put(300,15){\line(1,0){64}} \put(308,7){\line(0,1){32}}
\put(324,7){\line(0,1){16}} \multiput(300,15)(8,0){9}{\cc}
\multiput(308,7)(0,8){5}{\cc} \multiput(324,7)(0,16){2}{\cc}
\put(330,2){\tiny 12}

\put(380,15){\line(1,0){64}} \put(436,7){\line(0,1){32}}
\put(428,31){\line(1,0){16}} \multiput(380,15)(8,0){9}{\cc}
\multiput(436,7)(0,8){5}{\cc} \multiput(428,31)(16,0){2}{\cc}
\put(410,2){\tiny 13} \end{picture}\]

\[\begin{picture}(410,60) \put(0,30){\line(1,0){48}}
\put(8,6){\line(0,1){48}} \put(24,22){\line(0,1){16}}
\multiput(0,30)(8,0){7}{\cc} \multiput(8,6)(0,8){7}{\cc}
\multiput(24,22)(0,16){2}{\cc} \put(23,2){\tiny 14}

\put(60,15){\line(1,0){64}} \put(100,7){\line(0,1){32}}
\put(116,7){\line(0,1){16}} \multiput(60,15)(8,0){9}{\cc}
\multiput(100,7)(0,8){5}{\cc} \multiput(116,7)(0,16){2}{\cc}
\put(80,2){\tiny 15}

\put(140,15){\line(1,0){48}} \put(148,7){\line(0,1){16}}
\put(180,7){\line(0,1){48}} \multiput(140,15)(8,0){7}{\cc}
\multiput(148,7)(0,16){2}{\cc} \multiput(180,7)(0,8){7}{\cc}
\put(162,2){\tiny 16}

\put(200,15){\line(1,0){64}} \put(208,7){\line(0,1){16}}
\put(224,7){\line(0,1){32}} \multiput(200,15)(8,0){9}{\cc}
\multiput(224,7)(0,8){5}{\cc} \multiput(208,7)(0,16){2}{\cc}
\put(240,2){\tiny 17}

\put(280,15){\line(1,0){64}} \put(288,7){\line(0,1){32}}
\put(280,31){\line(1,0){16}} \multiput(280,15)(8,0){9}{\cc}
\multiput(288,7)(0,8){5}{\cc} \multiput(280,31)(16,0){2}{\cc}
\put(320,2){\tiny 18}

\put(360,15){\line(1,0){48}} \put(368,7){\line(0,1){48}}
\put(400,7){\line(0,1){16}} \multiput(360,15)(8,0){7}{\cc}
\multiput(368,7)(0,8){7}{\cc} \multiput(400,7)(0,16){2}{\cc}
\put(382,2){\tiny 19} \end{picture}\]

\[\begin{picture}(280,25) \put(0,15){\line(1,0){80}}
\put(24,7){\line(0,1){16}} \put(56,7){\line(0,1){16}}
\multiput(0,15)(8,0){11}{\cc} \multiput(24,7)(0,16){2}{\cc}
\multiput(56,7)(0,16){2}{\cc} \put(38,2){\tiny 20}

\put(100,15){\line(1,0){80}} \put(108,7){\line(0,1){16}}
\put(140,7){\line(0,1){16}} \multiput(100,15)(8,0){11}{\cc}
\multiput(108,7)(0,16){2}{\cc} \multiput(140,7)(0,16){2}{\cc}
\put(120,2){\tiny 21}

\put(200,15){\line(1,0){80}} \put(208,7){\line(0,1){16}}
\put(272,7){\line(0,1){16}} \multiput(200,15)(8,0){11}{\cc}
\multiput(208,7)(0,16){2}{\cc} \multiput(272,7)(0,16){2}{\cc}
\put(238,2){\tiny 22} \end{picture}\] and enumerate trees:
\[\begin{picture}(385,20) \put(0,8){\tiny $\alpha:$}
\put(15,2){\line(2,1){24}} \put(23,14){\line(2,-1){24}}
\put(31,10){\line(0,1){8}} \multiput(15,2)(8,4){4}{\cc}
\multiput(23,14)(8,-4){4}{\cc} \put(31,18){\cc}

\put(60,8){\tiny $\beta:$} \put(75,2){\line(2,1){24}}
\put(83,14){\line(2,-1){24}} \put(91,10){\line(0,-1){8}}
\multiput(75,2)(8,4){4}{\cc} \multiput(83,14)(8,-4){4}{\cc}
\put(91,2){\cc}

\put(120,8){\tiny $\gamma:$} \put(135,10){\line(1,0){40}}
\put(143,2){\line(0,1){16}} \multiput(135,10)(8,0){6}{\cc}
\multiput(143,2)(0,16){2}{\cc}

\put(190,8){\tiny $\delta:$} \put(205,10){\line(1,0){40}}
\put(221,2){\line(0,1){16}} \multiput(205,10)(8,0){6}{\cc}
\multiput(221,2)(0,16){2}{\cc}

\put(260,8){\tiny $\lambda:$} \put(275,10){\line(1,0){40}}
\put(283,10){\line(0,1){8}} \put(299,10){\line(0,-1){8}}
\multiput(275,10)(8,0){6}{\cc} \put(283,18){\cc} \put(299,2){\cc}

\put(330,8){\tiny $\mu:$} \put(345,10){\line(1,0){40}}
\put(353,10){\line(0,1){8}} \put(369,10){\line(0,1){8}}
\multiput(345,10)(8,0){6}{\cc} \put(353,18){\cc} \put(369,18){\cc}
\end{picture}\]
\[\begin{picture}(145,15) \put(0,3){\tiny $\nu:$}
\put(15,5){\line(1,0){48}} \put(39,5){\line(0,1){8}}
\multiput(15,5)(8,0){7}{\cc} \put(39,13){\cc}

\put(80,3){\tiny $\rho:$} \put(95,5){\line(1,0){48}}
\put(103,5){\line(0,1){8}} \multiput(95,5)(8,0){7}{\cc}
\put(103,13){\cc} \end{picture}\] As in the case of degree 6, we
for each degenerated cactus will: a) enumerate cacti to which
border it belongs; b) enumerate trees that can be obtained from it
by contraction.
$$\begin{array}{llllll} \hspace{2mm}1:\,\, A,L;\, \alpha,\lambda.&
\hspace{2mm}2:\,\, A,Q;\, \lambda,\mu.& \hspace{2mm}3:\,\, A,O;\,
\alpha,\mu. & \hspace{2mm}4:\,\, B,F;\,\beta,\mu. &
\hspace{2mm}5:\,\, B,S;\,\delta,\mu. & \hspace{2mm}6:\,\, B,C;\,
\beta,\delta.\\ \hspace{2mm}7:\,\, C,R;\,\delta,\lambda. &
\hspace{2mm}8:\,\, C,H;\,\beta,\lambda. & \hspace{2mm}9:\,\,
F,T;\, \beta,\delta. & 10:\,\,F,S;\,\mu,\rho. & 11:\,\,H,R;\,
\lambda,\rho. & 12:\,\,H,T;\,\beta,\rho.\\
13:\,\,J,P;\,\gamma,\nu. & 14:\,\, J,M;\, \alpha,\gamma. & 15:\,\,
J,O;\,\alpha,\nu. & 16:\,\,L,Q;\,\lambda,\nu. & 17:\,\, L,M;\,
\alpha,\nu. & 18:\,\,M,P;\,\gamma,\nu.\\ & 19:\,\,O,Q;\,\mu,\nu. &
20:\,\,P,P;\,\gamma,\gamma.& 21:\,\,R,S;\,\delta,\rho.&
22:\,\,T,T;\, \rho,\rho.& \end{array}$$ \begin{theor}
Combinatorial moduli space of Zolotarev polynomials of degree 7
with passport $(2,2,2)$ and the Galois $A_7$ has two components.
Each component is a projective map of 7 faces, 11 edges and 5
vertices.
\end{theor} \pn The first component: faces $A,J,L,M,O,P,Q$; edges
$1,2,3,13,14,15,16,17,18,19,20$; vertices
$\alpha,\gamma,\lambda,\mu,\nu$. The second component: faces
$B,C,F,H,R,S,T$; edges $4,5,6,7,8,9,10,11,12,21,22$; vertices
$\beta,\delta,\lambda,\mu,\rho$. In figure below we present maps
of both components.
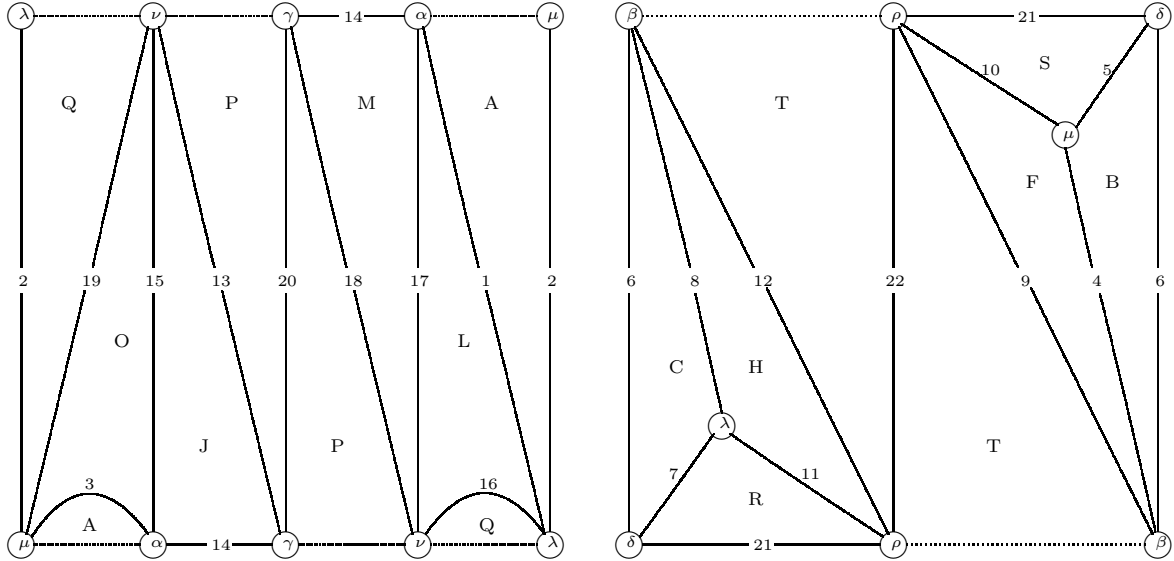
\[\begin{picture}(450,210) \multiput(5,5)(50,0){5}{\circle{10}}
\multiput(5,205)(50,0){5}{\circle{10}}
\multiput(5,10)(50,0){5}{\line(0,1){90}}
\multiput(5,110)(50,0){5}{\line(0,1){90}}
\multiput(4,4)(200,200){2}{\tiny $\mu$}
\multiput(4,204)(200,-200){2}{\tiny $\lambda$}
\multiput(54,204)(100,-200){2}{\tiny $\nu$}
\multiput(54,4)(100,200){2}{\tiny $\alpha$}
\multiput(104,4)(0,200){2}{\tiny $\gamma$}
\multiput(4,103)(200,0){2}{\tiny 2} \qbezier[40](10,5)(30,5)(50,5)
\qbezier[40](10,205)(30,205)(50,205)
\qbezier[40](160,5)(180,5)(200,5)
\qbezier[40](160,205)(180,205)(200,205)
\qbezier[40](110,5)(130,5)(150,5)
\qbezier[40](60,205)(80,205)(100,205)
\multiput(60,5)(25,0){2}{\line(1,0){15}}
\multiput(110,205)(25,0){2}{\line(1,0){15}}
\multiput(77,3)(50,200){2}{\tiny 14} \put(102,103){\tiny 20}
\put(52,103){\tiny 15} \put(152,103){\tiny 17}
\qbezier(103,9)(92,55)(81,101) \qbezier(57,201)(67,155)(78,110)
\qbezier(153,9)(142,55)(131,101)
\qbezier(107,201)(117,155)(128,110) \put(77,103){\tiny 13}
\put(127,103){\tiny 18} \put(82,170){\scriptsize P}
\put(122,40){\scriptsize P} \put(132,170){\scriptsize M}
\put(72,40){\scriptsize J} \qbezier(7,9)(18,55)(29,101)
\qbezier(53,201)(42,155)(32,110) \put(28,103){\tiny 19}
\put(40,80){\scriptsize O} \qbezier(203,9)(192,55)(181,101)
\qbezier(157,201)(167,155)(178,110) \put(179,103){\tiny 1}
\put(170,80){\scriptsize L} \put(180,170){\scriptsize A}
\put(20,170){\scriptsize Q} \qbezier(9,8)(30,40)(53,9)
\put(28,10){\scriptsize A} \put(29,26){\tiny 3}
\qbezier(158,9)(180,40)(203,9) \put(178,10){\scriptsize Q}
\put(178,26){\tiny 16}

\multiput(235,5)(100,0){3}{\circle{10}}
\multiput(235,205)(100,0){3}{\circle{10}}
\multiput(235,10)(100,0){3}{\line(0,1){90}}
\multiput(235,110)(100,0){3}{\line(0,1){90}}
\multiput(234,4)(200,200){2}{\tiny $\delta$}
\multiput(234,204)(200,-200){2}{\tiny $\beta$}
\multiput(334,4)(0,200){2}{\tiny $\rho$}
\multiput(234,103)(200,0){2}{\tiny 6}
\qbezier[50](240,205)(285,205)(330,205)
\qbezier[50](340,5)(385,5)(430,5) \put(332,103){\tiny 22}
\multiput(240,5)(50,0){2}{\line(1,0){40}}
\multiput(340,205)(50,0){2}{\line(1,0){40}}
\multiput(282,3)(100,200){2}{\tiny 21}
\qbezier(237,201)(260,155)(282,110)
\qbezier(287,102)(310,55)(333,9) \put(282,103){\tiny 12}
\qbezier(337,201)(360,155)(382,110)
\qbezier(387,102)(410,55)(433,9) \put(383,103){\tiny 9}
\put(290,170){\scriptsize T} \put(370,40){\scriptsize T}
\put(270,50){\circle{10}} \put(400,160){\circle{10}}
\put(269,49){\tiny $\lambda$} \put(399,159){\tiny $\mu$}
\qbezier(236,200)(247,155)(258,110)
\qbezier(260,100)(265,77)(270,55) \put(258,103){\tiny 8}
\qbezier(434,10)(424,55)(413,100)
\qbezier(400,155)(405,132)(410,110) \put(410,103){\tiny 4}
\qbezier(239,8)(253,27)(267,47) \qbezier(273,47)(302,27)(331,8)
\put(250,70){\scriptsize C} \put(280,70){\scriptsize H}
\put(280,20){\scriptsize R} \put(250,30){\tiny 7}
\put(300,30){\tiny 11} \qbezier(431,202)(417,182)(404,163)
\qbezier(338,202)(367,183)(397,164) \put(390,185){\scriptsize S}
\put(385,140){\scriptsize F} \put(415,140){\scriptsize B}
\put(414,183){\tiny 5} \put(368,183){\tiny 10}
\end{picture}\]
\begin{rem} Let $p=\int (x^2-1)(x^2+ax+b)(x^2+cx+d)\,dx$. The elimination
of variables $c$ and $d$ from the system
$$\begin{array}{l} p(1)=p(-1)\\ \text{the remainder $r$, where
$p=s\cdot (x^2+ax+b)+r$ has degree 0}\\ \text{the remainder g,
where $p=h\cdot(x^2+cx+d)+g$ has degree 0}\end{array}$$ gives us
the polynomial in $a$ and $b$ of degree 20. This polynomial
factorizes over $\mathbb{Q}$  in two factors of degrees 6 and 14,
respectively. The factor of degree 6 describes polynomials with
the Galois group $PSL_2(7)$. The factor of degree 14 factorizes
over $\mathbb{Q}(\sqrt{21})$ in two factors of degree 7. Each of
them defines a curve of genus 1. \end{rem}

\vspace{5mm}
\end{document}